\documentclass[aos,MSNbibl,nameyear,noautosecdot,dvips]{arximspdf}
\usepackage{graphicx}
\doi{10.1214/13-AOS1189} 
\volume{42}
\issue{1}
\pubyear{2014}
\firstpage{285}
\lastpage{323}

\makeatletter

\newcommand{\rright}{\right}
\newcommand{\lleft}{\left}
\newcommand{\rrvert}{\vert}
\newcommand{\llvert}{\vert}
\newtheorem{theorem}{Theorem}
\newtheorem{lemma}{Lemma}
\newtheorem{corollary}{Corollary}
\newproclaim{rem}{Remark}
\newcommand{\argmin}{\mathop{\mathrm{argmin}}}
\makeatother

\begin{document}
\begin{frontmatter}

\title{Adaptive piecewise polynomial estimation\\ via trend filtering\thanksref{T1}}
\runtitle{Trend filtering}

\pdftitle{Adaptive piecewise polynomial estimation via trend filtering}

\begin{aug}
\author[A]{\fnms{Ryan J.} \snm{Tibshirani}\corref{}\ead[label=e1]{ryantibs@cmu.edu}}
\runauthor{R. J. Tibshirani}
\affiliation{Carnegie Mellon University}
\address[A]{Department of Statistics\\
Carnegie Mellon University\\
Baker Hall\\
Pittsburgh, Pennsylvania 15213\\
USA\\
\printead{e1}} 
\end{aug}
\thankstext{T1}{Supported by NSF Grant DMS-13-09174.}

\received{\smonth{10} \syear{2013}}
\revised{\smonth{11} \syear{2013}}

%
\begin{abstract}
We study trend filtering, a recently proposed tool of
Kim et al. [\textit{SIAM Rev.} \textbf{51} (2009) 339--360]
for nonparametric regression. The trend filtering
estimate is defined as the minimizer of a penalized least squares
criterion, in which the penalty term sums the absolute $k$th order
discrete derivatives over the input points. Perhaps not surprisingly,
trend filtering estimates appear to have the structure of $k$th degree
spline functions, with adaptively chosen knot points (we say
``appear'' here as trend filtering estimates are not really
functions over continuous domains, and are only defined over
the discrete set of inputs). This brings
to mind comparisons to other nonparametric regression tools that also
produce adaptive splines; in particular, we compare trend filtering to
smoothing splines, which penalize the sum of squared derivatives
across input points, and to locally adaptive regression splines
[\textit{Ann. Statist.} \textbf{25} (1997) 387--413], which penalize the total variation of the
$k$th derivative. Empirically, we discover that
trend filtering estimates adapt to the local level of
smoothness much better than smoothing splines, and
further, they exhibit a remarkable similarity
to locally adaptive regression splines. We also provide
theoretical support for these empirical findings; most notably, we
prove that (with the right choice of tuning parameter)
the trend filtering estimate converges to the true
underlying function at the minimax rate for functions whose $k$th
derivative is of bounded variation. This is done via an asymptotic
pairing of trend filtering and locally adaptive regression splines,
which have already been shown to converge at the minimax rate
[\textit{Ann. Statist.} \textbf{25} (1997) 387--413].
At the core of this argument is a new result tying
together the fitted values of two lasso problems that share the same
outcome vector, but have different predictor matrices.
\end{abstract}

%

\pdfkeywords{62G08, 62G20, Trend filtering,
nonparametric regression,
smoothing splines,
locally adaptive regression splines,
minimax convergence rate,
lasso stability}

\end{frontmatter}\newpage

\setcounter{footnote}{1}
\section{\texorpdfstring{Introduction.}{Introduction}}\label{secintro}

Per the usual setup in nonparametric regression, we assume
that we have observations $y_1,\ldots, y_n \in\mathbb{R}$ from the model
%
\begin{equation}
\label{eqmodel} y_i = f_0(x_i) +
\varepsilon_i,\qquad i=1,\ldots, n,
\end{equation}
where $x_1, \ldots, x_n\in\mathbb{R}$ are input points, $f_0$ is the
underlying function to be estimated, and $\varepsilon_1,\ldots,
\varepsilon_n$ are independent errors. For the most part, we will
further assume that the inputs are evenly spaced
over the interval $[0,1]$, that is, $x_i=i/n$ for $i=1,\ldots, n$.
(However, this
assumption can be relaxed, as discussed in the supplementary
document [\citet{trendfilter-supp}].) The literature on
nonparametric regression is rich and diverse, and there are many
methods for estimating $f_0$ given \mbox{observations} from the model
(\ref{eqmodel}); some well-known examples include methods based on
local polynomials, kernels, splines, sieves and wavelets.

\begin{figure}

\includegraphics{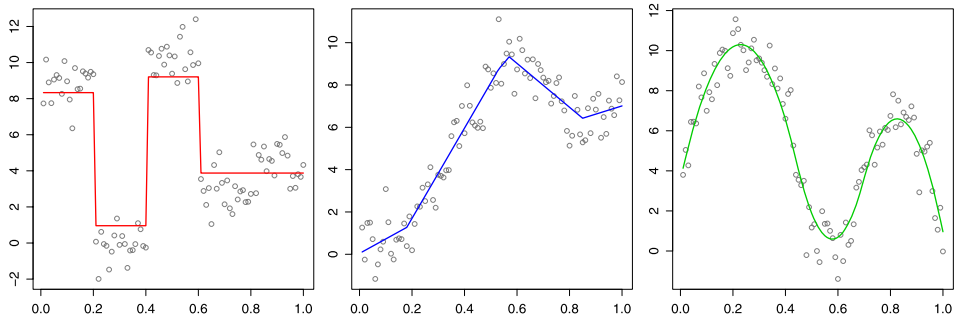}

\caption{Simple examples of trend filtering for constant,
linear, and quadratic orders ($k=0,1,2$, resp.), shown from left
to right. Although the trend
filtering estimates are only defined at the discrete inputs
$x_i=i/n$, $i=1,\ldots, n$, we use linear interpolation to extend the
estimates over $[0,1]$ for visualization purposes (this is the default
for all figures in this paper).}
\label{figtfsimp}
\end{figure}

This paper focuses on a relative newcomer in nonparametric
regression: trend filtering, proposed by \citet{l1tf}. For a
given integer $k \geq0$, the $k$th order trend filtering
estimate $\hat{\beta}=(\hat{\beta}_1,\ldots,\hat{\beta}_n)$ of
$(f_0(x_1), \ldots,
f_0(x_n))$ is defined by a penalized least squares optimization
problem,
%
\begin{equation}
\label{eqtf} \hat{\beta}= \argmin_{\beta\in\mathbb{R}^n} \frac{1}{2}\|y-\beta
\|_2^2 + \frac{n^k}{k!} \cdot\lambda \bigl\|D^{(k+1)}
\beta\bigr\|_1,
\end{equation}
where $\lambda\geq0$ is a tuning parameter, and
$D^{(k+1)}\in\mathbb{R}^{(n-k-1)\times n}$ is the discrete difference operator
of order $k+1$. (The constant factor $n^k/k!$ multiplying $\lambda$ is
unimportant, and can be absorbed into the tuning parameter
$\lambda$, but it will facilitate comparisons in future sections.)
When $k=0$, 
%
\begin{equation}
\label{eqd1} D^{(1)} = \lleft[ \matrix{ -1 & 1 & 0 & \cdots& 0 &
0
\vspace*{2pt}\cr
0 & -1 & 1 & \cdots& 0 & 0
\vspace*{2pt}\cr
\vdots& & & & &
\vspace*{2pt}\cr
0 & 0 & 0 & \cdots& -1 &
1}\rright] \in\mathbb{R}^{(n-1)\times n}
\end{equation}
and so $\|D^{(1)}\beta\|_1=\sum_{i=1}^{n-1} = | \beta_i-\beta_{i+1}|$.
Hence, the 0th order trend filtering
problem, which we will also call constant trend
filtering, is the same as one-dimensional total variation denoising
[\citet{tv}], or the one-dimensional fused lasso [\citet{fuse}]
(with pure fusion penalty, i.e., without an additional $\ell_1$
penalty on the coefficients themselves). In this case, $k=0$,
the components of the trend filtering estimate form a
piecewise constant structure, with\vspace*{1pt} points corresponding to the
nonzero entries of $D^{(1)}\hat{\beta}= (\hat{\beta}_2-\hat{\beta
}_1,\ldots,
\hat{\beta}_n-\hat{\beta}_{n-1})$. See Figure~\ref{figtfsimp} for
an example.

For $k \geq1$, the operator $D^{(k+1)} \in\mathbb{R}^{(n-k-1)\times
n}$ is
most easily-defined recursively, as in
%
\begin{equation}
\label{eqdk} D^{(k+1)}=D^{(1)}\cdot D^{(k)}.
\end{equation}
[Above, $D^{(1)}$ is the $(n-k-1)\times(n-k)$ version of the
first order discrete difference operator (\ref{eqd1}).]
In words, the definition (\ref{eqdk}) says that
the $(k+1)$st order difference operator is built up by evaluating
differences of differences, a total of $k+1$ times. Therefore, the
matrix $D^{(k+1)}$ can be thought of as the discrete analogy to the
$(k+1)$st order derivative operator, and the penalty term in
(\ref{eqtf}) penalizes the discrete $(k+1)$st derivative of the
vector $\beta\in\mathbb{R}^n$, that is, the changes in the discrete $k$th
derivative of $\beta$. Accordingly, one might expect the components
of the $k$th order trend filtering estimate to exhibit the
structure of a piecewise polynomial of order $k$, for example, for first order
trend filtering, the estimate would be piecewise linear, for second
order, it would be piecewise quadratic, etc.
Figure~\ref{figtfsimp} gives empirical evidence towards this claim.
Later, in Section~\ref{seccont}, we provide a more definitive
confirmation of this piecewise polynomial structure when we examine a
continuous-time representation for trend filtering.

%

It is straightforward to check that
\begin{eqnarray*}
D^{(2)} &=& \lleft[ \matrix{ 1 & -2 & 1 & 0 & \cdots& 0
\vspace*{2pt}\cr
0 & 1 & -2 &
1 &\cdots& 0
\vspace*{2pt}\cr
0 & 0 & 1 & -2 & \cdots& 0
\vspace*{2pt}\cr
\vdots& & & & &} \rright],
\\
D^{(3)} &=& \lleft[ %
\matrix{ -1 & 3 & -3 & 1 & \cdots& 0
\vspace*{2pt}\cr
0
& -1 & 3 & -3 & \cdots& 0
\vspace*{2pt}\cr
0 & 0 & -1 & 3 & \cdots& 0
\vspace*{2pt}\cr
\vdots& & & & &}
\rright]
\end{eqnarray*}
and in general, the nonzero elements in each row of $D^{(k)}$ are
given by the \mbox{$(k+1)$st} row of Pascal's triangle, but with alternating
signs. A more explicit (but also more complicated-looking)
expression for the $k$th order trend filtering problem is therefore
\[
\hat{\beta}= \argmin_{\beta\in\mathbb{R}^n} \frac{1}{2}\sum
_{i=1}^n (y_i-\beta_i)^2
+ \frac{n^k}{k!} \cdot \lambda \sum
_{i=1}^{n-k-1} \Biggl\llvert \sum
_{j=i}^{i+k+1} (-1)^{j-i} \pmatrix{k+1 \cr j-i}
\beta_j \Biggr\rrvert.
\]
The penalty term above sums over successive linear combinations
of $k+2$ adjacent coefficients, that is, the discrete difference operator
$D^{(k+1)}$ is a banded matrix with bandwidth $k+2$.

\subsection{\texorpdfstring{The generalized lasso and related properties.}{The generalized lasso and related properties}}\label{secglprops}

\begin{figure}

\includegraphics{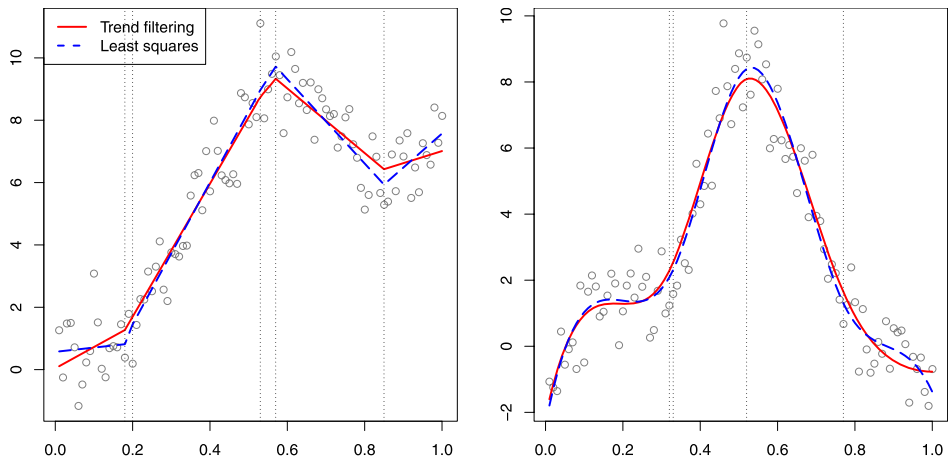}

\caption{Examples of the shrinkage effect for
linear trend filtering ($k=1$, left panel) and for cubic trend filtering
($k=3$, right panel). In each panel, the solid red line is the trend
filtering estimate (at a particular value of $\lambda$), and the
dashed blue line is the regression spline estimate of the same order
and with the same knots, with the vertical lines marking the
locations of the knots. The trend filtering estimates on the left and
right have shrunken 1st and 3rd derivatives, respectively, compared
to their regression spline counterparts.}
\label{figtfshrink}
\end{figure}

For any order $k\geq0$, the trend filtering estimate $\hat{\beta}$ is
uniquely defined, because the criterion in (\ref{eqtf}) is strictly
convex. Furthermore, the trend filtering criterion is of
generalized lasso form with an identity predictor matrix $X=I$ (this
is called the signal approximator case)
and a specific choice of penalty matrix $D=D^{(k+1)}$. Some properties
of the trend filtering estimate therefore follow from known results
on the generalized lasso [\citeauthor{genlasso} (\citeyear{genlasso,lassodf2})], for example, an
exact representation of the trend filtering estimate in terms of
its active set and signs,
and also, a formula for its degrees of freedom:
%
\begin{equation}
\label{eqdf} \mathrm{df}(\hat{\beta}) = \mathbb{E}[\mbox{number of knots in $
\hat{\beta}$}] + k + 1,
\end{equation}
where the number of knots in $\hat{\beta}$ is interpreted to mean
the number of nonzero entries in $D^{(k+1)}\hat{\beta}$ (the
basis and continuous-time representations of trend filtering, in
Sections~\ref{sectflasso} and \ref{seccont},
provide a justification for this interpretation).
To repeat some of the discussion in \citeauthor{genlasso} (\citeyear{genlasso,lassodf2}), the result in (\ref{eqdf}) may seem
somewhat remarkable, as a fixed-knot $k$th
degree regression spline with $d$ knots also has $d+k+1$ degrees of
freedom---and trend filtering does not employ fixed knots, but
rather, chooses them adaptively. So, why does trend filtering not have
a larger degrees of freedom? At a high level, the answer lies in the
shrinkage due to the $\ell_1$ penalty in (\ref{eqtf}): the nonzero
entries of $D^{(k+1)}\hat{\beta}$ are shrunken toward zero,
compared to the same quantity for the corresponding
equality-constrained least squares estimate.
In other words, within each interval defined by the (adaptively
chosen) knots, trend filtering fits a $k$th degree polynomial whose
$k$th derivative is shrunken toward its $k$th derivatives in neighboring
intervals, when compared to a $k$th degree regression spline with the
same knots. Figure~\ref{figtfshrink} gives a demonstration of this
phenomenon for $k=1$ and $k=3$.

%

In terms of algorithms,
the fact that the discrete difference operator $D^{(k+1)}$ is banded
is of great advantage for solving the generalized lasso problem in
(\ref{eqtf}). \citet{l1tf} describe a primal--dual interior point
method for solving (\ref{eqtf}) at a fixed value of $\lambda$,
wherein each iteration essentially reduces to
solving a linear system in $D^{(k+1)}(D^{(k+1)})^T$, costing $O(n)$
operations. In the worst case, this algorithm requires $O(n^{1/2})$
iterations, so its complexity is $O(n^{3/2})$.\footnote{It should
be noted that hidden in the $O(\cdot)$ notation here is a factor
depending on the prespecified
error tolerance $\varepsilon$, namely, a term of the form
$\log(1/\varepsilon)$. We emphasize here that the primal--dual
interior point method is a different type of algorithm than
the path algorithm, in the sense that the latter returns an exact
solution up to computer precision, whereas the former returns
an $\varepsilon$-suboptimal solution, as measured by the difference in its
achieved criterion value and the optimal criterion value. Essentially,
all general purpose convex optimization techniques (that are
applicable to the trend filtering problem) fall into the same class
as the primal--dual interior point method, that is, they return
$\varepsilon$-suboptimal solutions; only specialized techniques like the
path algorithm can deliver exact solutions.} [\citet{l1tf} focus
mainly on linear
trend filtering, the case $k=1$, but their arguments carry over to the
general case as well.] On the other hand, instead of solving
(\ref{eqtf}) at a fixed~$\lambda$, \citet{genlasso} describe a
path algorithm to solve (\ref{eqtf}) over all values of $\lambda\in
[0,\infty)$, that is, to compute the entire solution path
$\hat{\beta}=\hat{\beta}(\lambda)$ over $\lambda$. This path is
piecewise linear as a function of $\lambda$ (not to be
confused with the estimate itself at any
fixed $\lambda$, which has a piecewise polynomial
structure over the input points $x_1,\ldots, x_n$).
Again, the bandedness of $D^{(k+1)}$ is key here for efficient
computations, and \citet{genlasso-alg} describe an
implementation of the path algorithm in which computing the estimate
at each successive critical point in the path requires $O(n)$ operations.

Software for both of these algorithms is freely available
online. For the primal--dual interior point method, see
\url{http://stanford.edu/\textasciitilde boyd/l1_tf}, which provides Matlab and C
implementations (these only cover the linear trend
filtering case, but can be extended to the general polynomial
case); for the path algorithm, see the function \texttt{trendfilter} in
the R package \texttt{genlasso}, available on the CRAN repository.

\subsection{\texorpdfstring{Summary of our results.}{Summary of our results}}\label{secqanda}

\begin{figure}

\includegraphics{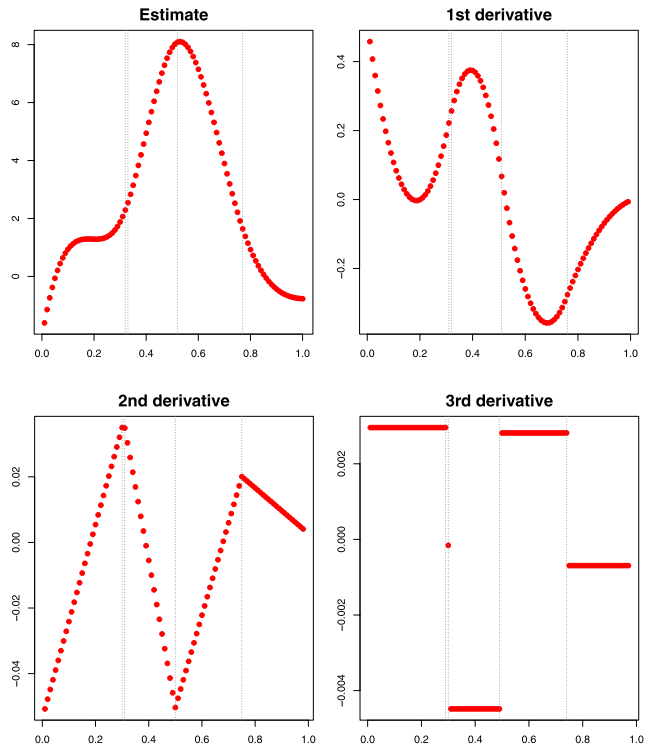}

\caption{The leftmost panel shows the same
cubic trend filtering estimate as in Figure~\protect\ref{figtfshrink}
(but here we do not use linear interpolation to
emphasize the discrete nature of the estimate). The components
of this estimate appear to be the evaluations of a continuous
piecewise polynomial function. Moreover, its discrete 1st and
2nd derivatives (given by multiplying by $D^{(1)}$ and $D^{(2)}$,
resp.) also appear to be continuous, and its
discrete third derivative (from multiplication by $D^{(3)}$) is
piecewise constant within each interval. Hence, we might believe
that such a trend filtering estimate actually represents the
evaluations of a 3rd degree spline over the inputs
$x_i=i/n$, $i=1,\ldots, n$. We address this idea in Section~\protect\ref{seccont}.}\label{figtfspline}
\end{figure}

Little is known about trend filtering---mainly, the results
due to its generalized lasso form, for example, the degrees
of freedom result (\ref{eqdf}) discussed in the previous
section---and much is unknown. Examining the trend
filtering fits in Figures~\ref{figtfsimp} and \ref{figtfshrink}, it
appears that the estimates not only have the structure
of piecewise polynomials, they furthermore have the structure of
splines:
these are piecewise polynomial functions that have continuous
derivatives of all orders lower than the leading one [i.e., a $k$th
degree spline is a $k$th degree piecewise polynomial with continuous
$0$th through $(k-1)$st derivatives at its knots].
Figure~\ref{figtfspline} plots an example cubic trend filtering
estimate, along with its discrete 1st, 2nd and 3rd derivatives
(given by multiplication by $D^{(1)}$, $D^{(2)}$, and $D^{(3)}$,
resp.). Sure enough, the lower order discrete derivatives
appear ``continuous'' across the knots, but what does this really mean
for such discrete sequences? Does trend filtering have an analogous
continuous-time representation, and if so, are the estimated functions
really splines?

%

Besides these questions, one may also wonder about the
performance of trend filtering estimates compared to other methods.
Empirical examples (like those in Section~\ref{secss})
show that trend filtering estimates achieve a significantly higher
degree of local adaptivity than smoothing splines, which are arguably
the standard tool for adaptive spline estimation. Other examples (like
those in Section~\ref{seclrs}) show that trend filtering estimates
display a comparable level of adaptivity to locally adaptive regression
splines, another well-known technique for adaptive spline estimation,
proposed by
\citet{locadapt} on the basis of being more locally adaptive (as
their name would suggest). Examples are certainly
encouraging, but a solely empirical conclusion here would be
unsatisfactory---fixing as a metric the squared error loss in
estimating the true function $f_0$, averaged over the input points,
can we say more definitively that trend filtering estimates actually
outperform smoothing splines, and perform as well as locally adaptive
regression splines?

We investigate the questions discussed above in this paper.
To summarize our results, we find that:
\begin{itemize}
\item for $k=0,1$ (constant or linear orders), the
continuous-time analogues of trend filtering estimates are indeed
$k$th degree splines; moreover, they are exactly the same as $k$th
order locally adaptive regression splines;
\item for $k \geq2$ (quadratic or higher orders), the continuous-time
versions of trend filtering estimates are not quite splines, but
piecewise polynomial functions that are ``close to'' splines (with
small discontinuities in lower order derivatives at the
knots); hence, they are not the same as $k$th order locally adaptive
regression splines;
\item for any $k$, if the $k$th derivative of true function $f_0$ is
of bounded variation, then the $k$th order trend filtering
estimate converges to $f_0$ (in terms of squared error loss) at the
minimax rate; this rate is achieved by locally adaptive regression
splines [\citet{locadapt}], but not by smoothing splines nor any other
estimate linear in $y$ [\citet{minimaxwave}].
\end{itemize}
We note that, although trend filtering and locally adaptive regression
splines are formally different estimators for $k\geq2$, they are
practically indistinguishable by eye in most examples.
Such a degree of similarity, in finite samples, goes beyond what we
are able to show theoretically. However, we do prove that
trend filtering estimates and locally adaptive regression spline
estimates converge to each other asymptotically (Theorem
\ref{thmtflrs}). The argument here boils down to a bound on
the difference in the fitted values of two lasso problems
that have the same outcome vector, but different predictor matrices
(because both trend filtering and locally adaptive regression
splines can be represented as lasso problems, see Section~\ref{seclrs}). To the best of our knowledge, this general bound is
a new result (see the supplementary document [\citet{trendfilter-supp}]).
Further, we use this asymptotic pairing between trend
filtering and locally adaptive regression splines to prove the minimax
convergence rate for trend filtering (Corollary \ref{cortfrate}).
The idea is simple: trend filtering and locally adaptive regression
splines converge to each other at the minimax rate, locally adaptive
regression splines converge to the true function at the minimax rate
[\citet{locadapt}], and hence so does trend filtering.

\subsection{\texorpdfstring{Why trend filtering?}{Why trend filtering?}}\label{secwhytf}


Trend filtering estimates, we argue, enjoy the favorable
theoretical performance of locally adaptive regression splines; but
now it is fair to ask: why would we ever use trend
filtering, over, say, the latter estimator? The main reason is
that trend filtering estimates are much easier to compute, due to
the bandedness of the discrete derivative operators, as explained
previously. The computations for locally adaptive regression splines,
meanwhile, cannot exploit such sparsity or structure, and are
considerably slower. To be more concrete, the
primal--dual interior point method described in Section~\ref
{secglprops} above can handle problems of size on the order of
$n={}$1,000,000 points (and the path algorithm, on the order of
$n={}$100,000 points), but even for $n={}$10,000 points, the computations
for locally adaptive regression splines are prohibitively slow.
We discuss this in Section~\ref{seclrs}.

Of course, the nonparametric regression toolbox is highly-developed
and already offers plenty of good methods.
We do not presume that trend filtering should be regarded as the
preferred method in every nonparametric regression problem, but simply
that it represents a useful contribution to the toolbox, being both fast
and locally adaptive, that is, balancing the strengths of smoothing
splines and locally adaptive regression splines. This manuscript
mainly focuses on the comparison to the aforementioned estimators
because they, too, like trend filtering,
fit piecewise polynomials functions and they are widely used.
Though we do not compare wavelets or smoothing splines with a
spatially variable tuning parameter in as much detail, we consider
them in Section~\ref{secastro} in an analysis of astrophysics data.
It should be mentioned that for trend filtering to become a truly
all-purpose nonparametric regression tool, it must be able to handle
arbitrary input points $x_1,\ldots, x_n$ (not just evenly spaced
inputs). We give an extension to this case in the supplementary
document [\citet{trendfilter-supp}]. Our analysis of trend filtering
with arbitrary inputs\vadjust{\goodbreak} shows promising computational and theoretical
properties, but still, a few questions remain unanswered. This will
be the topic of future work.

As a separate point, another distinguishing feature of trend
filtering is that it falls into what is called
the \textit{analysis} framework with respect to its problem formulation,
whereas locally adaptive regression splines, smoothing splines, and
most others fall into the \textit{synthesis} framework. Synthesis
and analysis are two terms used in signal processing that
describe different approaches for defining an
estimator with certain desired characteristics. In the synthesis
approach, one builds up the estimate constructively from
a set of characteristic elements or atoms; in the analysis approach,
the strategy is instead to define the estimate deconstructively, via
an operator that penalizes undesirable or uncharacterisic
behavior. Depending on the situation,
it can be more natural to implement the former
rather than the latter, or vice versa, and hence both are important.
We discuss the importance of the analysis framework in the context of
nonparametric regression estimators in Section~\ref{secsva}, where we
define extensions of trend filtering that would be difficult to
construct from the synthesis perspective, for example, a sparse variant of
trend filtering.

Here is an outline for the rest of this article (though we have
discussed its contents throughout the \hyperref[secintro]{Introduction}, we list them here
in proper order). In Sections~\ref{secss} and
\ref{seclrs}, we compare trend filtering to smoothing splines and
locally adaptive regression splines, respectively. We give data
examples that show trend filtering estimates are more locally
adaptive than smoothing splines, and that trend filtering and locally
adaptive regression splines are remarkably similar, at any common
value of their tuning parameters. We
also discuss the differing computational requirements for these
methods. In Section~\ref{seclrs}, we show that both locally adaptive
regression splines and trend filtering can be posed as lasso problems,
with identical predictor matrices when $k=0$ or 1, and
with similar but slightly different predictor matrices when $k\geq2$.
This allows us to conclude that trend filtering and locally adaptive
regression splines are exactly the same for constant or linear orders,
but not for quadratic or higher orders. Section~\ref{seccont}
develops a continuous-time representation for the trend filtering
problem, which reveals that (continuous-time) trend filtering
estimates are always $k$th order piecewise polynomials,
but for $k \geq2$, are not $k$th order splines. In Section~\ref
{secrates}, we derive the minimax convergence rate of trend
filtering estimates, under the assumption that the $k$th derivative of
the true function has bounded total variation. We do this by
bounding the difference between trend filtering estimates and locally
adaptive regression splines, and invoking the fact that the latter are
already known to converge at the minimax rate [\citet{locadapt}].
We also study convergence rates for a true function with growing
total variation. In Section~\ref{secastro}, we consider an
astrophysics data example, and compare the performance of several
commonly used nonparametric regression tools.
Section~\ref{secextdisc} presents ideas for future work:
multivariate trend filtering, sparse trend filtering, and the
synthesis versus analysis perspectives.
Essentially all proofs, and the discussion of trend filtering
for arbitrary inputs, are deferred until the supplementary
document [\citet{trendfilter-supp}].

\section{\texorpdfstring{Comparison to smoothing splines.}{Comparison to smoothing splines}}\label{secss}

Smoothing splines are a popular tool in nonparametric regression, and
have been extensively studied in terms of both computations and theory
[some well-known references are \citet{deboorsplines},
\citet{wahbasplines}, \citet{greensilver}]. Given input
points $x_1,\ldots, x_n \in[0,1]$, which we assume are unique, and
observations $y_1,\ldots, y_n \in\mathbb{R}$, the $k$th order
smoothing spline
estimate is defined as
%
\begin{equation}
\label{eqss1} \hat{f}= \argmin_{f \in\mathcal{W}_{(k+1)/2}} \sum
_{i=1}^n \bigl(y_i-f(x_i) \bigr)^2 + \lambda\int_0^1
\bigl(f^{((k+1)/2)}(t) \bigr)^2 \,dt,
\end{equation}
where $f^{((k+1)/2)}(t)$ is the derivative of $f$ of order
$(k+1)/2$, $\lambda\geq0$ is a tuning parameter, and the domain of
minimization here is the Sobolev space
\begin{eqnarray*}
\mathcal{W}_{(k+1)/2} &=& \biggl\{ f\dvtx [0,1]\rightarrow\mathbb{R}\dvtx  f \mbox{ is
} (k+1)/2\mbox{ times differentiable and }
\\
&&\hspace*{126pt} \int_0^1 \bigl(f^{((k+1)/2)}(t) \bigr)^2 \,dt < \infty \biggr\}.
\end{eqnarray*}
Unlike trend filtering, smoothing splines are only defined for an odd
polynomial order $k$. In practice, it seems that the case $k=3$ (i.e.,
cubic smoothing splines) is by far the most common case considered.
In the next section, we draw a high-level comparison between
smoothing splines and trend filtering, by writing the smoothing spline
minimization problem (\ref{eqss1}) in finite-dimensional
form. Following this, we make empirical comparisons, and then discuss
computational efficiency.

\subsection{\texorpdfstring{Generalized ridge regression and Reinsch form.}{Generalized ridge regression and Reinsch form}}

Remarkably, it can be shown that the infinite-dimensional problem in
(\ref{eqss1}) is has a unique minimizer, which is a $k$th degree
natural spline with knots at the input points $x_1,\ldots, x_n$ [see,
e.g., \citet{wahbasplines}, \citet{greensilver},
\citet{esl}]. Recall that a $k$th degree natural spline is a simply
a $k$th degree spline that reduces to a polynomial of degree $(k-1)/2$
before the first knot and after the last knot; it is easy to check the
set of natural splines of degree $k$, with knots at $x_1,\ldots, x_n$,
is spanned by precisely $n$ basis functions. Hence, to solve
(\ref{eqss1}), we can solve for the coefficients $\theta\in\mathbb
{R}^n$ in
this basis expansion:
%
\begin{equation}
\label{eqss2} \hat{\theta}= \argmin_{\theta\in\mathbb{R}^n} \|y-N\theta
\|_2^2 + \lambda\theta^T \Omega\theta,
\end{equation}
where $N \in\mathbb{R}^{n\times n}$ contains the evaluations of
$k$th degree natural spline basis functions over the knots
$x_1,\ldots, x_n$, and $\Omega\in\mathbb{R}^{n\times n}$ contains\vadjust{\goodbreak} the
integrated products of their $((k+1)/2)$nd derivatives; that is, if
$\eta_1,\ldots,\eta_n$ denotes a collection of basis functions for
the set
of $k$th degrees natural splines with knots at $x_1,\ldots, x_n$, then
%
\begin{eqnarray}\label{eqssbp}
N_{ij} = \eta_j(x_i) \quad\mbox{and}\quad \Omega_{ij} = \int_0^1
\eta_i^{((k+1)/2)}(t) \cdot\eta_j^{((k+1)/2)}(t)
\,dt
\nonumber\\[-9pt]\\[-9pt]
\eqntext{\mbox{for all } i,j=1,\ldots, n.}
\end{eqnarray}
The problem in (\ref{eqss2}) is a generalized ridge regression, and
from its solution $\hat{\theta}$, the function $\hat{f}$ in
(\ref{eqss1}) is simply given at the input points $x_1,\ldots, x_n$
by
\[
\bigl(\hat{f}(x_1),\ldots,\hat{f}(x_n) \bigr) = N\hat{
\theta}.
\]
More generally, the smoothing spline estimate $\hat{f}$ at an
arbitrary input $x\in[0,1]$ is given by
\[
\hat{f}(x) = \sum_{j=1}^n \hat{
\theta}_j \eta_j(x).
\]

To compare the smoothing spline problem, as expressed in
(\ref{eqss2}), with trend filtering, it helps to rewrite the
smoothing spline fitted values as follows:
%
\begin{eqnarray}\label{eqreinsch}
N\hat{\theta}&=& N\bigl(N^T N + \lambda\Omega
\bigr)^{-1} N^T y\nonumber
\\
&=& N \bigl(N^T \bigl(I + \lambda N^{-T} \Omega
N^{-1}\bigr) N \bigr)^{-1} N^T y
\\
&=& (I+\lambda K)^{-1} y,\nonumber
\end{eqnarray}
where $K=N^{-T}\Omega N^{-1}$. The expression in (\ref{eqreinsch}) is
called the \textit{Reinsch} form for the
fitted values. From this expression, we can view $\hat{u}=N\hat
{\theta}$
as the solution of the minimization problem
%
\begin{equation}
\label{eqss3} \hat{u} = \argmin_{u\in\mathbb{R}^n} \|y-u\|_2^2
+ \lambda u^T K u,
\end{equation}
which is of similar form to the trend filtering problem in (\ref{eqtf}),
but here the $\ell_1$ penalty $\|D^{(k+1)}\beta\|_1$
is replaced by the quadratic penalty $u^T K u = \|K^{1/2}
u\|_2^2$. How do these two penalties compare? First,
the penalty matrix $K^{1/2}$ used by smoothing splines is similar in
nature to the discrete derivative operators [we know from its
continuous-time analog in (\ref{eqss1}) that the term
$\|K^{1/2}u\|_2^2$ penalizes wiggliness in something like the
$((k+1)/2)$nd derivative of $u$] but is still strictly different.
For example, for $k=3$
(cubic smoothing splines) and input points $x_i=i/n$, $i=1,\ldots, n$,
it can be shown [\citet{greenyandell}] that\vspace*{1pt} the smoothing spline penalty
is $\|K^{1/2} u\|_2^2 = \|C^{-1/2} D^{(2)} u\|_2^2/n^3$ where
$D^{(2)}$ is the second order discrete derivative operator, and
$C\in\mathbb{R}^{n\times n}$ is a tridiagonal matrix, with diagonal
elements equal to $2/3$ and off-diagonal elements equal to $1/6$.

A second and more important difference is that smoothing splines
utilize a (squared) $\ell_2$ penalty, while trend filtering uses an
$\ell_1$ penalty. Analogous to the usual comparisons between
ridge regression and the lasso, the former penalty shrinks the
components of $K^{1/2} \hat{u}$, but does not set any of the
components to zero unless $\lambda=\infty$ (in which case all
components are zero), whereas the latter penalty shrinks and also
adaptively sets components of $D^{(k+1)}\hat{\beta}$ to zero. One might
imagine, recalling that $K^{1/2}$ and $D^{(k+1)}$ both act in a
sense as derivative operators, that trend
filtering estimates therefore exhibit a finer degree of local
adaptivity than do smoothing splines.
This idea is supported by the examples in the next section, which
show that trend filtering estimates outperform smoothing splines (when
both are optimally tuned) in estimating functions with spatially
inhomogeneous smoothness. The idea is also supported by our theory in
Section~\ref{secrates}, where we prove that trend filtering estimates
have a better rate of convergence than smoothing splines
(in fact, they achieve the optimal rate) over a broad class of
underlying functions.\looseness=1

\subsection{\texorpdfstring{Empirical comparisons.}{Empirical comparisons}}\label{secssempir}

We compare trend filtering and smoothing spline estimates on
simulated data. We fix $k=3$ (i.e., we compare cubic trend
filtering versus cubic smoothness splines), because the
\texttt{smooth.spline} function in the R programming language provides
a fast implementation for smoothing splines in this case. Generally
speaking, smoothing
splines and trend filtering provide similar estimates when the
underlying function $f_0$ has spatially homogeneous smoothness, or
to put it simply, is either entirely smooth or entirely wiggly
throughout its domain. Hence, to illustrate the difference between
the two estimators, we consider two examples of functions that
display varying levels of smoothness at different spatial locations.

%
\begin{figure}

\includegraphics{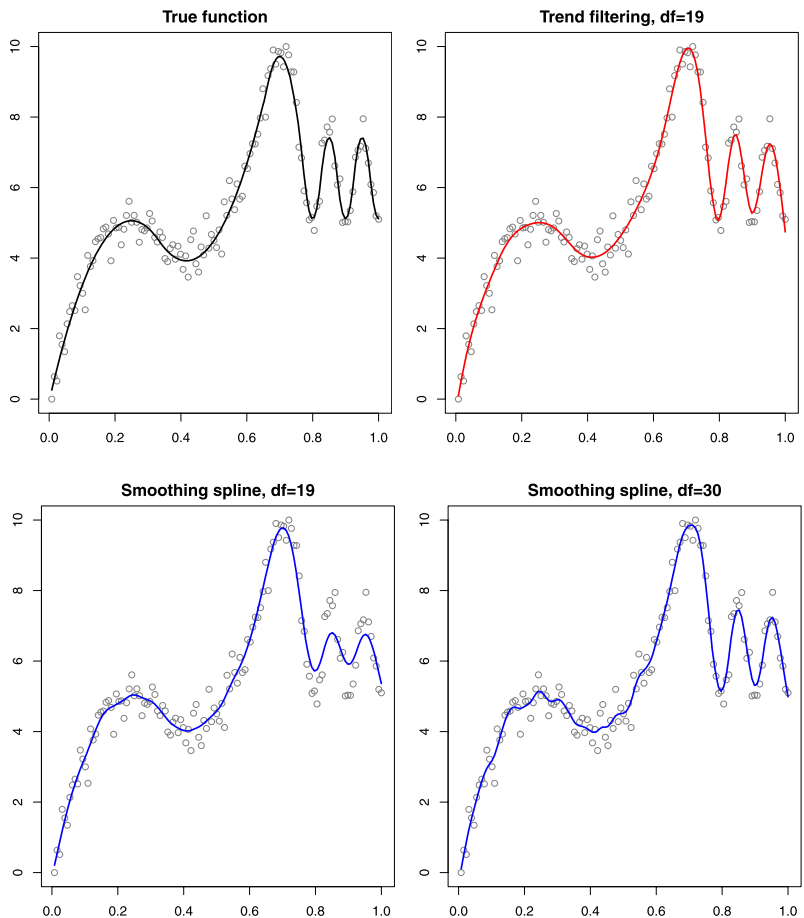}

\caption{An example with $n=128$ observations drawn from a model
where the underlying function has variable spatial smoothness, as
shown in the top left panel. The cubic trend filtering estimate with
19 degrees of freedom, shown in the top right panel, picks up
the appropriate level of smoothness at different spatial locations:
smooth at the left-hand side of the domain, and wiggly at the right-hand
side. When also allowed 19 degrees of freedom, the cubic smoothing
spline estimate in the bottom left panel
grossly underestimates the signal on the right-hand side of the domain. The
bottom right panel shows the smooth spline estimate with 30 degrees of
freedom, tuned so that it displays the appropriate level of adaptivity
on the right-hand side; but now, it is overly adaptive on the left-hand side.}
\label{figsmoothwiggly}
\end{figure}

Our first example, which we call the ``hills'' example, considers
a piecewise cubic function $f_0$ over $[0,1]$, whose knots
are spaced farther apart on the left-hand side of the domain, but bunched
closer together on the right-hand side. As a
result, $f_0(x)$ is smooth for $x$ between $0$ and about $0.8$, but
then abruptly becomes more wiggly---see the top left panel of
Figure~\ref{figsmoothwiggly}. We drew $n=128$ noisy observations from $f_0$
over the evenly spaced inputs $x_i=i/n$, $i=1,\ldots, n$ (with
independent, normal noise), and fit a trend filtering estimate, tuned
to have 19 degrees of freedom, as shown in the top right
panel.\footnote{To be precise, this is an unbiased estimate of its
degrees of freedom; see (\ref{eqdf}) in Section~\ref{secglprops}.}
We can see here that the estimate adapts to the appropriate levels of
smoothness at both the left and right sides of the domain. But this is
not true of the smoothing spline estimate with 19 degrees of freedom,
displayed in the bottom left panel: the estimate is considerably
oversmoothed on
the right-hand side. As we increase the allowed flexibility, the smoothing
spline estimate is able to fit the small hills on the right, with a
total of 30 degrees of freedom; however, this causes undersmoothing on
the left-hand side, as shown in the bottom right panel.

%
\begin{figure}

\includegraphics{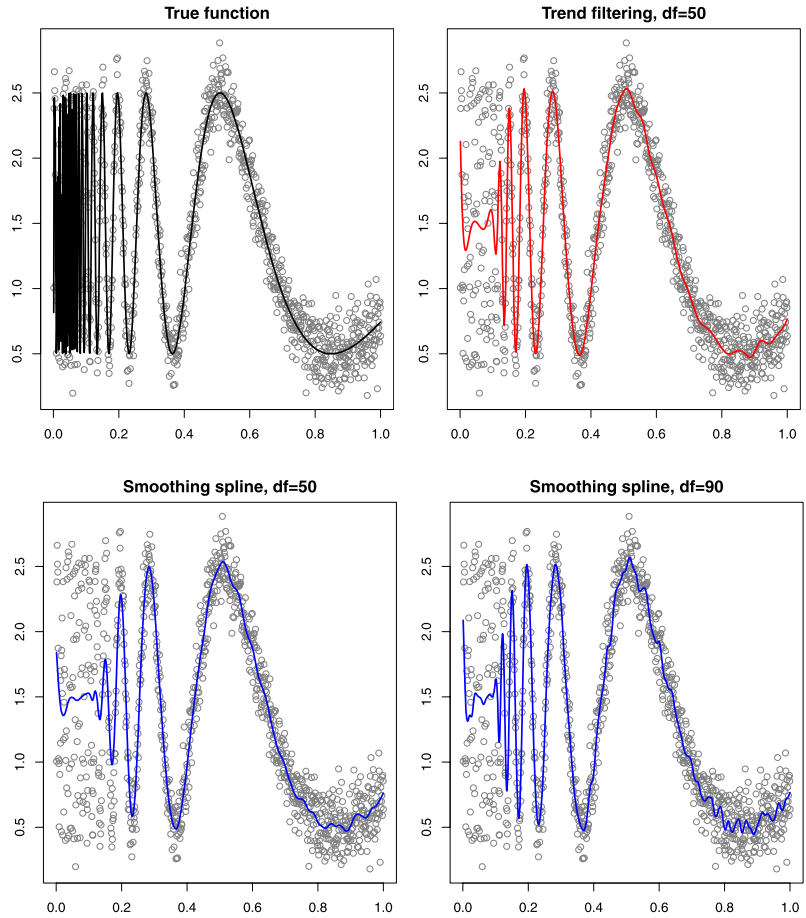}

\caption{An example with $n=1000$ noisy observations of the
Doppler function, $f(x)= \sin(4/x)+1.5$, drawn in the top left
panel. The top right and bottom left panels show the cubic trend
filtering and smoothing spline estimates, each with 50 degrees of
freedom; the former captures approximately 4 cycles of the
Doppler function, and the latter only 3. If we nearly double the
model complexity, namely, we use 90 degrees of freedom, then
the smoothing spline estimate is finally able to capture 4 cycles,
but the estimate now becomes very jagged on the right-hand side of
the plot.}\label{figdoppler}
\end{figure}

For our second example, we take $f_0$ to be the ``Doppler'' function
[as considered in, e.g., \citet{sure},
\citet{locadapt}]. Figure~\ref{figdoppler}, clockwise from the
top left, displays the Doppler function and corresponding $n=1000$
noisy observations, the trend filtering estimate with 50 degrees of
freedom, the smoothing spline estimate with 50 degrees of freedom,
and the smoothing spline estimate with 90 degrees of freedom. The same
story, as in the hills example, holds here: trend filtering adapts to
the local level of smoothness better than smoothing splines,
which have trouble with the rapidly increasing frequency of the
Doppler function (as $x$ decreases).

%
\begin{figure}

\includegraphics{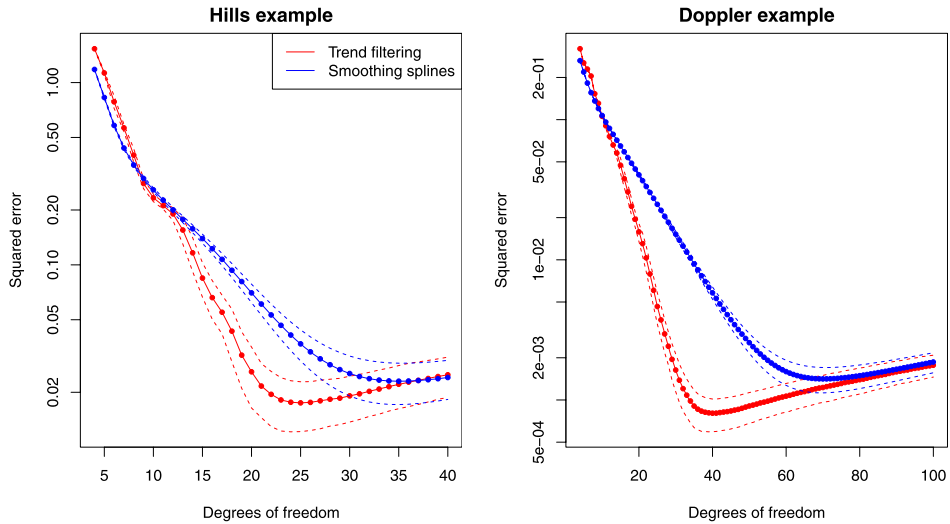}

\caption{Shown is the squared error loss in
predicting the true function $f_0$, averaged over the
input points, for the hills data example on the left, and the
Doppler example on the right. In each setup, trend filtering and
smoothing spline estimators were fit over a range of degrees of
freedom values; the red curves display the loss for trend filtering,
and the blue curves for smoothing splines. The results were averaged
over 50 simulated data sets, and the standard deviations are denoted
by dotted lines. In these examples, trend filtering has a generally
better predictive accuracy than smoothing splines, especially for
models of low to intermediate complexity (degrees of freedom).}\label{figerrs}
\end{figure}

In Figure~\ref{figerrs}, we display the input-averaged squared
error losses\footnote{For the Doppler data example, we actually
average the
squared error loss only over inputs $x_i \geq0.175$, because for
$x_i<0.175$, the true Doppler function $f_0$ is of such high frequency
that neither trend filtering nor smoothing splines are able to do a
decent job of fitting it.}
\[
\frac{1}{n} \sum_{i=1}^n \bigl(
\hat{\beta}_i-f_0(x_i)\bigr)^2
\quad\mbox{and}\quad \frac{1}{n} \sum_{i=1}^n
\bigl(\hat{f}(x_i)-f_0(x_i)
\bigr)^2
\]
for the trend filtering and smoothing spline estimates
$\hat{\beta}$ and $\hat{f}$, respectively, on the hills and
Doppler examples. We considered a wide range of model complexities
indexed by degrees of freedom, and averaged the results
over 50 simulated data sets for each setup (the dotted lines show plus
or minus one standard deviations). Aside from the visual evidence
given in Figures~\ref{figsmoothwiggly} and \ref{figdoppler},
Figure~\ref{figerrs} shows that from the perspective of squared
error loss, trend filtering outperforms smoothing splines in
estimating underlying functions with variable spatial smoothness.
As mentioned previously, we will prove in Section~\ref{secrates} that
for a large class of underlying functions $f_0$, trend filtering estimates
have a sharper convergence rate than smoothing splines.

\subsection{\texorpdfstring{Computational considerations.}{Computational considerations}}\label{secsscomp}

Recall that the smoothing spline fitted values are given by
%
\begin{equation}
\label{eqssfit} N\hat{\theta}= N\bigl(N^T N + \lambda\Omega
\bigr)^{-1} N^T y,
\end{equation}
where $N \in\mathbb{R}^{n\times n}$ contains the
evaluations of basis functions $\eta_1,\ldots,\eta_n$ for the subspace
of $k$th degree natural splines with knots at the inputs, and
$\Omega\in\mathbb{R}^{n\times n}$ contains their integrated products
of their
$((k+1)/2)$nd order derivatives, as in~(\ref{eqssbp}). Depending on
exactly which basis we choose, computation of (\ref{eqssfit}) can be
fast or slow; by choosing the $B$-spline basis functions, which have
local support, the matrix $N^T N + \lambda\Omega$ is banded, and
so the smoothing spline fitted values can be computed in $O(n)$
operations [e.g., see \citet{deboorsplines}]. In practice, these
computations are extremely fast.


By comparison, \citet{l1tf} suggest a primal--dual interior point
method, as mentioned in Section~\ref{secglprops}, that computes the
trend filtering estimate (at any fixed value of the tuning parameter
$\lambda$) by iteratively solving a sequence of banded linear
systems, rather than just\vspace*{1pt} a single one. Theoretically, the
worst-case number of iterations scales as $O(n^{1/2})$, but the
authors report that in practice the number of iterations needed is
only a few tens, almost independent of the problem size~$n$.
Hence, trend filtering computations with the
primal--dual path interior point method are
slower than those for smoothing splines, but not by a huge margin.

To compute the trend filtering estimates for
the examples in the previous section, we actually used the dual path
algorithm of \citet{genlasso}, which was also discussed in
Section~\ref{secglprops}. Instead of solving the trend
filtering problem at a fixed value of $\lambda$, this
algorithm constructs the solution path as $\lambda$ varies
from~$\infty$ to $0$. Essentially, it does so by stepping
through a sequence of estimates, where each step either adds one knot to
or deletes one knot from the fitted piecewise polynomial
structure. The computations at each step amount to solving two banded
linear systems, and hence require $O(n)$ operations; the overall
efficiency depends on how many steps along the path are
needed before the estimates of interest have been reached (at which
point the path algorithm can be terminated early). But because
knots can be both added and deleted to the fitted piecewise polynomial
structure at each step, the algorithm can take much more than
$k$ steps to reach an estimate with $k$ knots. Consider the
Doppler data example, in the last section, with $n=1000$ points: the
path algorithm used nearly 4000 steps to compute the trend filtering
estimate with 46 knots (50 degrees of freedom) shown in the upper
right panel of Figure~\ref{figdoppler}. This took approximately 28
seconds on a standard desktop computer, compared to the smoothing
spline estimates shown in the bottom left and\vadjust{\goodbreak} right panels of
Figure~\ref{figdoppler}, which took about 0.005 seconds each. We
reiterate that in this period of time, the path algorithm for
trend filtering computed a total of 4000 estimates, versus a single
estimate computed by the smoothing spline solver. (A quick
calculation, $28/4000 = 0.007$, shows that the
time per estimate here is comparable.) For the hills data
set in the last section, where $n=128$, the dual path
algorithm constructed the entire path of trend filtering
estimates (consisting of 548 steps) in less than 3 seconds; both
smoothing spline estimates took under 0.005 seconds each.

\section{\texorpdfstring{Comparison to locally adaptive regression splines.}{Comparison to locally adaptive regression splines}}\label{seclrs}

Locally adaptive regression splines are an alternative
to smoothing splines, proposed by \citet{locadapt}. They are more
computationally intensive than smoothing splines but have better
adaptivity properties (as their name would suggest). Let $x_1,\ldots,
x_n\in[0,1]$ denote the inputs, assumed unique and ordered as in
$x_1<x_2<\cdots<x_n$,
and $y_1,\ldots, y_n \in\mathbb{R}$ denote the observations. For the $k$th
order locally adaptive regression spline estimate, where $k\geq0$ is
a given arbitrary integer (not necessarily odd, as is required for
smoothing splines), we start by defining the knot superset
%
\begin{equation}
\label{eqt} T = \cases{ \{x_{k/2+2},\ldots, x_{n-k/2}\}, &\quad
if $k$ is even,
\vspace*{2pt}\cr
\{x_{(k+1)/2+1},\ldots, x_{n-(k+1)/2}\}, &\quad if $k$
is odd.}
\end{equation}
This is essentially just the set of inputs $\{x_1,\ldots, x_n\}$, but
with points near the left and right boundaries removed. We then define
the $k$th order locally adaptive regression spline estimate as
%
\begin{equation}
\label{eqlrs1} \hat{f}= \argmin_{f \in\mathcal{G}_k} \frac{1}{2}\sum
_{i=1}^n \bigl(y_i-f(x_i)
\bigr)^2 + \lambda\cdot \mathrm{TV}\bigl(f^{(k)}\bigr),
\end{equation}
where $f^{(k)}$ is now the $k$th weak derivative
of $f$, $\lambda\geq0$ is a tuning parameter,
$\mathrm{TV}(\cdot)$ denotes the total variation operator, and
$\mathcal{G}_k$ is the set
%
\begin{equation}\label{eqgk}
\mathcal{G}_k = \bigl\{ f\dvtx [0,1]\rightarrow\mathbb{R}\dvtx
f \mbox{ is $k$th degree spline with knots contained in $T$} \bigr\}.\hspace*{-35pt}
\end{equation}
Recall that for a function $f\dvtx  [0,1] \rightarrow\mathbb{R}$, its total
variation is defined as
\[
\mathrm{TV}(f) = \sup \Biggl\{\sum_{i=1}^p\bigl|f(z_{i+1})-f(z_i)\bigr|\dvtx  \mbox{$z_1 < \cdots< z_p$ is a partition of $[0,1]$}
\Biggr\}
\]
and this reduces to $\mathrm{TV}(f) = \int_0^1 |f'(t)| \,dt$ if $f$ is
(strongly) differentiable.

Next, we briefly address the difference between our definition of
locally adaptive regression splines in (\ref{eqlrs1}) and the
original definition found in \citet{locadapt}; this discussion
can be skipped without interrupting the flow of ideas.
After this, we rewrite problem (\ref{eqlrs1}) in terms of
the coefficients of $f$ with respect to a basis for the
finite-dimensional set~$\mathcal{G}_k$. For an arbitrary choice
of basis, this new problem is of generalized lasso form, and in
particular, if we choose the truncated power series as our basis for
$\mathcal{G}_k$, it simply becomes a lasso problem. We will
see that trend filtering, too, can be represented as a
lasso problem, which allows for a more direct comparison between
the two estimators.

\subsection{\texorpdfstring{Unrestricted locally adaptive regression splines.}{Unrestricted locally adaptive regression splines}}

For readers familiar with the work of \citet{locadapt},
it may be helpful to explain the difference between our definition of
locally adaptive regression splines and theirs: these authors
define the locally adaptive regression spline estimate as
%
\begin{equation}
\label{equlrs} \hat{f}\in\argmin_{f \in\mathcal{F}_k} \frac{1}{2}\sum
_{i=1}^n \bigl(y_i-f(x_i)
\bigr)^2 + \lambda\cdot \mathrm{TV}\bigl(f^{(k)}\bigr),
\end{equation}
where $\mathcal{F}_k$ is the set
\[
\mathcal{F}_k = \bigl\{ f\dvtx [0,1]\rightarrow\mathbb{R}\dvtx  \mbox {$f$ is
$k$ times weakly differentiable and } \mathrm{TV}\bigl(f^{(k)}\bigr) <
\infty \bigr\}.
\]
[The element notation in (\ref{equlrs}) emphasizes the fact that the
minimizer is not generally unique.] We call (\ref{equlrs}) the \textit{unrestricted} locally adaptive regression spline problem, in reference
to its minimization domain compared to that of (\ref{eqlrs1}).
\citet{locadapt} prove that the minimum in this
unrestricted problem is always achieved by a $k$th degree
spline, and that this spline has knots contained in $T$ if $k=0$ or
1, but could have knots outside of $T$ (and in fact, outside of the
input set $\{x_1,\ldots, x_n\}$) if $k \geq2$. In other words, the
solution in (\ref{eqlrs1}) is always a solution in (\ref{equlrs})
when $k=0$ or 1, but this need not be true when $k \geq2$; in the
latter case, even though there exists a $k$th degree spline that
minimizes (\ref{equlrs}), its knots could occur at noninput 
points.\looseness=-1

The unrestricted locally adaptive regression estimate
(\ref{equlrs}) is the main object of theoretical
study in \citet{locadapt}, but practically speaking, this
estimate is difficult to compute when $k\geq2$, because the possible
knot locations are generally not easy to determine [see also
\citet{plregpath}]. On the other hand, the restricted estimate
as defined in (\ref{eqlrs1}) is more computationally
tractable. Fortunately, \citet{locadapt} show that essentially
all of their theoretical results for the unrestricted estimate
also apply to the restricted estimate, as long as the input points
$x_1,\ldots, x_n$ are not spaced too far apart. In particular,
for evenly spaced inputs, $x_i=i/n$, $i=1,\ldots, n$, the convergence
rates of the unrestricted and restricted estimates are the same.
We therefore focus on the restricted problem (\ref{eqlrs1}) in the
current paper, and mention the unrestricted problem (\ref{equlrs})
purely out of interest. For example, to anticipate results to come, we
will show in Lemma \ref{lemsamediff} of Section~\ref{sectflasso}
that the trend filtering estimate (\ref{eqtf}) for $k=0$ or
1 is equal to the locally adaptive regression spline
estimate (\ref{eqlrs1}) (i.e., they match at the input points
$x_1,\ldots, x_n$); hence, from what we discussed above, it
also solves the unrestricted problem in (\ref{equlrs}).

\subsection{\texorpdfstring{Generalized lasso and lasso form.}{Generalized lasso and lasso form}}

We note that the knot set $T$ in (\ref{eqt}) has $n-k-1$
elements, so $\mathcal{G}_k$ is spanned by $n$ basis
functions, call them $g_1,\ldots, g_n$.
Since each $g_j$, $j=1,\ldots, n$ is a $k$th degree spline with knots
in $T$, we know that its $k$th weak derivative is piecewise constant
and (say) right-continuous, with jump points contained in $T$;
therefore, writing $t_0=0$ and $T=\{t_1,\ldots, t_{n-k-1}\}$, we have
\[
\mathrm{TV}\bigl(g_j^{(k)}\bigr) = \sum
_{i=1}^{n-k-1} \bigl| g_j^{(k)}(t_i)-g_j^{(k)}(t_{i-1})
\bigr|.
\]
Similarly, any linear combination of $g_1,\ldots, g_n$ satisfies
\[
\mathrm{TV} \Biggl( \Biggl(\sum_{j=1}^n
\theta_j g_j \Biggr)^{(k)} \Biggr) = \sum
_{i=1}^{n-k-1} \Biggl| \sum
_{j=1}^n \bigl(g_j^{(k)}(t_i)-g_j^{(k)}(t_{i-1})
\bigr)\cdot\theta_j \Biggr|.
\]
Hence, we can reexpress (\ref{eqlrs1}) in terms of the
coefficients $\theta\in\mathbb{R}^n$ in its basis expansion with
respect to
$g_1,\ldots, g_n$,
%
\begin{equation}
\label{eqlrs2} \hat{\theta}= \argmin_{\theta\in\mathbb{R}^n} \frac{1}{2}\| y-G
\theta\|_2^2 + \lambda \|C\theta\|_1,
\end{equation}
where $G \in\mathbb{R}^{n\times n}$ contains the evaluations of
$g_1,\ldots,
g_n$ over the inputs $x_1,\ldots, x_n$, and
$C \in\mathbb{R}^{(n-k-1)\times n}$ contains the differences in their $k$th
derivatives across the knots, that is,
%
\begin{eqnarray}
\label{eqg} G_{ij} &=& g_j(x_i)\qquad
\mbox{for } i,j=1,\ldots, n,
\\
\label{eqc} C_{ij} &=& g_j^{(k)}(t_i)-g_j^{(k)}(t_{i-1})
\qquad\mbox{for } i=1,\ldots, n-k-1, j=1,\ldots, n.
\end{eqnarray}
Given the solution $\hat{\theta}$ in (\ref{eqlrs2}), we can recover the
locally adaptive regression spline estimate $\hat{f}$ in (\ref{eqlrs1})
over the input points by
\[
\bigl(\hat{f}(x_1),\ldots,\hat{f}(x_n) \bigr) = G\hat{
\theta},
\]
or, at an arbitrary point $x\in[0,1]$ by
\[
\hat{f}(x) = \sum_{j=1}^n \hat{
\theta}_j g_j(x).
\]
The problem (\ref{eqlrs2}) is a generalized lasso problem,
with predictor matrix $G$ and penalty matrix $C$; by
taking $g_1,\ldots, g_n$ to be the truncated power basis, we can turn
(a block of) $C$ into the identity, and hence
(\ref{eqlrs2}) into a lasso problem.\looseness=-1

%
\begin{lemma}
\label{lemlassog}
Let $T=\{t_1,\ldots, t_{n-k-1}\}$ denote the set defined in
(\ref{eqt}), and let $g_1,\ldots, g_n$ denote the $k$th order
truncated power basis with knots in $T$,
%
\begin{eqnarray}\label{eqgbasis}
g_1(x)&=&1,\qquad g_2(x)=x, \ldots, g_{k+1}(x) = x^k,
\nonumber\\[-8pt]\\[-8pt]
g_{k+1+j}(x) &=& (x-t_j)^k \cdot1\{x \geq t_j\},\qquad j=1,\ldots, n-k-1.\nonumber
\end{eqnarray}
(For the case $k=0$, we interpret $0^0=1$.)
Then the locally adaptive regression spline
problem (\ref{eqlrs1}) is equivalent to the following lasso problem:
%
\begin{equation}
\label{eqlassog} \hat{\theta}= \argmin_{\theta\in\mathbb{R}^n} \frac{1}{2}\|y-G\theta
\|_2^2 + \lambda \sum
_{j=k+2}^n |\theta_j|,
\end{equation}
in that $\hat{f}(x) = \sum_{j=1}^n \hat{\theta}_j g_j(x)$ for $x\in[0,1]$.
Here, $G\in\mathbb{R}^{n\times n}$ is the basis matrix as in (\ref{eqg}).
\end{lemma}

Lemma \ref{lemlassog} follows from the fact that
for the truncated power basis, the penalty matrix $C$ in (\ref{eqc})
satisfies $C_{i,i+k+1}=1$ for $i=1,\ldots, n-k-1$, and $C_{ij}=0$
otherwise.
It is worth noting that \citet{lassoknots} investigate a lasso
problem similar to (\ref{eqlassog}) for the purposes of knot
selection in regression splines.

Note that (\ref{eqlassog}) is of somewhat nonstandard form for a
lasso problem, because the $\ell_1$ penalty is only taken
over the last $n-k-1$ components of $\theta$. We will see next that
the trend filtering problem in (\ref{eqtf}) can also be written in
lasso form (again with the $\ell_1$ penalty summing over the last
$n-k-1$ coefficients), and we will compare these two formulations.
First, however, it is helpful to express the knot superset
$T$ in (\ref{eqt}) and the basis matrix $G$ in (\ref{eqg}) in a more
explicit form, for evenly spaced input points $x_i=i/n$, $i=1,\ldots,
n$ (this being the underlying assumption for trend filtering). These
become
%
\begin{equation}
\label{eqtt}\quad  T = \cases{ \bigl((k+2)/2 + i\bigr)/n, &\quad for $i=1,\ldots,
n-k-1$, if $k$ is even,
\vspace*{2pt}\cr
\bigl((k+1)/2 + i\bigr)/n, &\quad for $i=1,\ldots,
n-k-1$, if $k$ is odd}
\end{equation}
and
%
\begin{equation}
\label{eqgg}  G_{ij} = \cases{
\cases{ 0, &\quad for $i < j$,
\vspace*{2pt}\cr
1, &
\quad for $i \geq j$,} &\quad if $k=0$,
\vspace*{5pt}\cr
\cases{ i^{j-1}/n^{j-1},
\vspace*{2pt}\cr
    \hspace*{32.5pt}\mbox{for $i=1,\ldots, n$, $j=1,\ldots, k+1$,}
\vspace*{2pt}\cr
0,\qquad\mbox{for $i \leq j-k/2$, $j \geq k+2$,}
\vspace*{2pt}\cr
(i-j+k/2)^k/n^k,
\vspace*{2pt}\cr
    \hspace*{32.5pt}\mbox{for $i > j-k/2$, $j \geq k+2$,}} &\quad if $k>0$ is even,
\vspace*{5pt}\cr
\cases{ i^{j-1}/n^{j-1},
\vspace*{2pt}\cr
\hspace*{32.7pt} \mbox{for $i=1,\ldots, n$, $j=1,\ldots, k+1$,}
\vspace*{2pt}\cr
0,\qquad\mbox{for $i \leq j-(k+1)/2$, $j \geq k+2$,}
\vspace*{2pt}\cr
\bigl(i-j+(k+1)/2\bigr)^k/n^k,
\vspace*{2pt}\cr
\hspace*{32.7pt} \mbox{for $i > j-(k+1)/2$, $j \geq k+2$,}} &\quad if $k>0$ is odd.}\hspace*{-35pt}
\end{equation}
(It is not really important to separate the definition of $G$ for
$k=0$ from that for $k>0$, $k$ even; this is only done to make
transparent the structure of $G$.)

\subsection{\texorpdfstring{Trend filtering in lasso form.}{Trend filtering in lasso form}}\label{sectflasso}

We can transform the trend filtering problem in (\ref{eqtf}) into
lasso form, just like the representation for locally adaptive
regression splines in (\ref{eqlassog}).

%
\begin{lemma}
\label{lemtflasso}
The trend filtering problem in (\ref{eqtf}) is equivalent to the
lasso problem
%
\begin{equation}
\label{eqlassoh} \hat{\alpha}= \argmin_{\alpha\in\mathbb{R}^n} \frac{1}{2}\|y-H\alpha
\|_2^2 + \lambda \sum_{j=k+2}^n |\alpha_j|,
\end{equation}
in that the solutions satisfy $\hat{\beta}= H\hat{\alpha}$.
Here, the predictor matrix $H\in\mathbb{R}^{n\times n}$ is given by
%
\begin{equation}
\label{eqhh} H_{ij} = \cases{ i^{j-1}/n^{j-1}, &
\quad for $i=1,\ldots, n$, $j=1,\ldots, k+1$,
\vspace*{2pt}\cr
0, &\quad for $i \leq j-1$, $j \geq
k+2$,
\vspace*{2pt}\cr
\sigma_{i-j+1}^{(k)} \cdot k!/n^k, &\quad
for $i > j-1$, $j \geq k+2$,}
\end{equation}
where we define $\sigma_i^{(0)}=1$ for all $i$ and
\[
\sigma_i^{(k)} = \sum_{j=1}^i
\sigma_j^{(k-1)}\qquad \mbox{for } k=1,2,3,\ldots,
\]
that is, $\sigma_i^{(k)}$ is the $k$th order cumulative sum of
$(1,1,\ldots,1)\in\mathbb{R}^i$.
\end{lemma}

The proof of this lemma basically inverts the $(k+1)$st
order discrete difference operator $D^{(k+1)}$; see the supplementary
document [\citet{trendfilter-supp}]. We remark that this result, in the
special case of $k=1$, can be found in \citet{l1tf}.

It is not hard to check that in the case $k=0$ or 1, the definitions
of $G$ in (\ref{eqgg}) and $H$ in (\ref{eqhh}) coincide,
which means that the locally adaptive regression spline and trend
filtering problems (\ref{eqlassog}) and (\ref{eqlassoh}) are the same.
But when $k \geq2$, we have $G\neq H$, and hence the problems are
different.

%
\begin{lemma}
\label{lemsamediff}
Consider evenly spaced inputs $x_i=i/n$, $i=1,\ldots, n$,
and the basis matrices $G,H$ defined in (\ref{eqgg}),
(\ref{eqhh}). If $k=0$ or 1, then $G=H$, so the lasso
representations for locally adaptive regression splines and trend
filtering, (\ref{eqlassog})~and~(\ref{eqlassoh}), are the same.
Therefore, their solutions are the same, or in other words,
\[
\hat{\beta}_i = \hat{f}(x_i) \qquad\mbox{for } i=1,
\ldots, n,
\]
where $\hat{\beta}$ and $\hat{f}$ are the solutions of the original trend
filtering and locally adaptive regression spline problems,
(\ref{eqtf}) and (\ref{eqlrs1}),
at any fixed common value of the tuning parameter $\lambda$.

If $k \geq2$, however, then $G \neq H$, so the problems
(\ref{eqlassog}) and (\ref{eqlassoh}) are different, and hence the
trend filtering and locally adpative regression spline estimators are
generically different.
\end{lemma}

See the supplement for the proof [\citet{trendfilter-supp}]. Though the
trend filtering and locally adaptive regression
spline estimates are formally different for polynomial
orders $k \geq2$, they are practically very similar (at all common
values of
$\lambda$). We give examples of this next, and then
compare the computational requirements for the two methods.

%
\begin{figure}[b]

\includegraphics{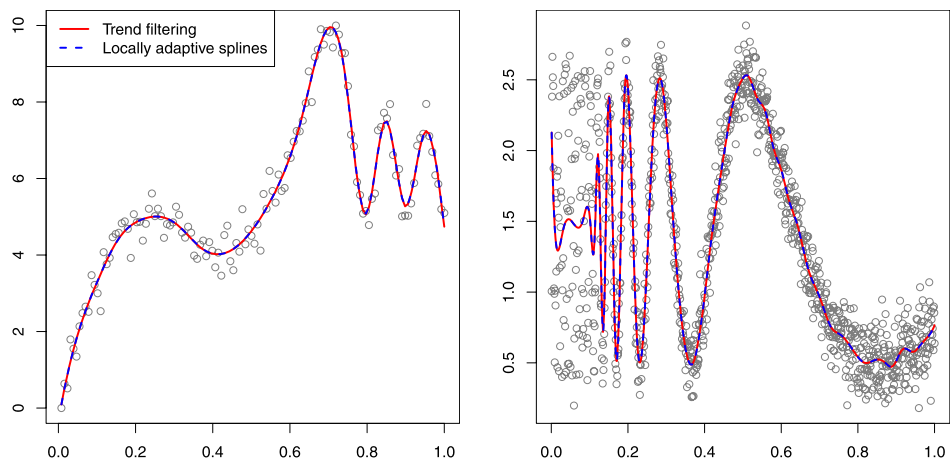}

\caption{Trend filtering and locally adaptive regression spline
estimates, using the same values of the tuning parameter $\lambda$, for
the hills and Doppler data examples considered in
Section~\protect\ref{secssempir}. The trend filtering estimates are
drawn as
solid red lines, and locally adaptive regression splines as dotted
blue lines; in both examples, the two estimates are basically
indistinguishable by eye.}\label{figlrs}
\end{figure}

\subsection{\texorpdfstring{Empirical comparisons.}{Empirical comparisons}}

We revisit the hills and Doppler examples of Section~\ref{secssempir}. Figure~\ref{figlrs} displays, for $k=3$ (cubic
order), the trend filtering and locally adaptive regression spline
estimates at matching values of the tuning parameter $\lambda$.
The estimates are visually identical in both examples (but not
numerically identical---the average squared difference between the
estimates across the input points is around $10^{-5}$ for
the hills example, and $10^{-7}$ for the Doppler example). This remains
true for a wide range of common tuning parameter values, and only for
very small values of $\lambda$ do slight differences between the two
estimators become noticeable.

Nothing is special about the choice $k=3$ here or about these data
sets in particular: as far as we can tell, the same phenomenon occurs
for any polynomial order~$k$, and any set of observations.
This extreme similarity between the two estimators, holding in finite
sample and across essentially all common tuning \mbox{parameter} values, is
beyond what we show theoretically in Section~\ref{secrates}. In this
section, we prove that for tuning parameters
of a certain order, the two estimators converge asymptotically at a
fast rate. Sharper statements are a topic for future work.

\subsection{\texorpdfstring{Computational considerations.}{Computational considerations}}
\label{seclrscomp}

Both the locally adaptive regression spline and trend filtering
problems can be represented as lasso problems with dense, square
predictor matrices, as in (\ref{eqlassog}) and (\ref{eqlassoh}).
For trend filtering, we do this purely for analytical reasons, and
computationally it is much more efficient to work from its original
representation in (\ref{eqtf}), where the penalty operator
$D^{(k+1)}$ is sparse and banded. As discussed in Sections~\ref{secglprops} and \ref{secsscomp}, two efficient algorithms for
trend filtering are the primal--dual interior point method of
\citet{l1tf} and the dual path algorithm of
\citet{genlasso}; the former computes the trend filtering
estimate at a fixed value of $\lambda$, in $O(n^{3/2})$ worst-case
complexity [the authors claim that the practical complexity is
closer to $O(n)$]; the latter computes the entire solution path
over $\lambda$, with each critical point in this piecewise linear path
requiring $O(n)$ operations.

For locally adaptive regression splines,
on the other hand, there is not a better computational alternative
than solving the lasso problem in (\ref{eqlassog}). One can
check that the inverse of the truncated power basis matrix $G$ is
dense, so if we converted~(\ref{eqlassog}) to generalized lasso form
[to match the form of trend filtering in (\ref{eqtf})], then it would
have a dense penalty matrix. And if we were to choose, for example, the
$B$-spline basis over the truncated power basis to parameterize
$\mathcal{G}_k(T)$ in (\ref{eqlrs1}), then although the basis matrix $G$
would be sparse and banded, the resulting penalty matrix $C$ in
(\ref{eqlrs2}) would be dense. In other words, to compute the
locally adaptive regression spline estimate, we are more or less
stuck with solving the lasso problem in (\ref{eqlassog}), where $G$
is the dense predictor matrix in (\ref{eqgg}). This task is
manageable for small or moderately sized problems, but for large
problems, dealing with the $n\times n$ dense matrix $G$, and even
holding it in memory, becomes burdensome.

To compute the locally adaptive regression spline estimates for the
examples in the last section, we
solved the lasso problem in (\ref{eqlassog}) using the LARS algorithm
for the lasso path [\citet{lars}], as implemented by the
\texttt{lars} R package.\footnote{To fit the problem in (\ref{eqlassog})
into standard lasso form, that is, a form in which the $\ell_1$ penalty
is taken over the entire coefficient
vector, we can solve directly for the first $k+1$ coefficients (in
terms of the last $n-k-1$ coefficients), simply by linear regression.}
This particular algorithm was chosen for the sake of a fair
comparison to the dual path algorithm used for trend filtering.
For the Doppler data example with $n=1000$ points, the LARS algorithm
computed the locally adaptive regression spline estimate (shown in the
right panel of Figure~\ref{figlrs}) in a comparable amount of time to
that taken by the dual path algorithm for trend filtering---in fact, it
was faster, at approximately 16 versus 28 seconds on a standard
desktop computer. The real issue, however,
is scalability. For $n=1000$ points, each of these algorithms required
about 4000 steps to compute their respective estimates; for $n={}$10,000
noisy observations from the Doppler curve, the dual path algorithm
completed 4000 steps in a little under 2.5 minutes, whereas the LARS
algorithm completed 4000 steps in 1 hour. Furthermore, for problem
sizes $n$ somewhat larger than $n={}$10,000, just fitting the $n\times
n$ basis matrix $G$ used by the LARS algorithm into memory becomes an
issue.

\section{\texorpdfstring{Continuous-time representation.}{Continuous-time representation}}\label{seccont}

Section~\ref{sectflasso} showed that the trend filtering minimization
problem (\ref{eqtf}) can be expressed in lasso form
(\ref{eqlassoh}), with a predictor matrix $H$ as in (\ref{eqhh}). The
question we consider is now: is there a set of
basis functions whose evaluations over the inputs $x_1,\ldots,
x_n$ give this matrix $H$? Our next lemma answers this question
affirmatively.

%
\begin{lemma}
\label{lemtfh}
Given inputs $x_1<\cdots<x_n$,
consider the functions $h_1,\ldots, h_n$ defined as
%
\begin{eqnarray}
\label{eqhbasis}
h_1(x) &=& 1,\qquad h_2(x) = x, \ldots, h_{k+1}(x) = x^k,
\nonumber\\[-8pt]\\[-8pt]
h_{k+1+j}(x) &=& \prod_{\ell=1}^k (x-x_{j+\ell}) \cdot1\{x \geq x_{j+k}\},\qquad j=1,\ldots, n-k-1.\nonumber
\end{eqnarray}
If the input points are evenly
spaced over $[0,1]$, $x_i=i/n$ for $i=1,\ldots, n$, then the trend
filtering basis matrix $H$ in (\ref{eqhh}) is generated by evaluating
the functions $h_1,\ldots, h_n$ over $x_1,\ldots, x_n$, that is,
%
\begin{equation}
\label{eqh} H_{ij} = h_j(x_i),\qquad i,j=1,
\ldots, n.
\end{equation}
\end{lemma}

The proof is given in the supplementary document [\citet{trendfilter-supp}].
As a result of the lemma, we can alternatively express the trend
filtering basis matrix $H$ in (\ref{eqhh}) as
%
\begin{equation}
\label{eqtfh}
H_{ij} = \cases{
i^{j-1}/n^{j-1},
\vspace*{2pt}\cr
    \hspace*{32.5pt}\mbox{for $i=1,\ldots, n$, $j=1,\ldots, k+1$,}
\vspace*{2pt}\cr
0,\qquad\mbox{for $i \leq j-1$, $j \geq k+2$,}
\vspace*{2pt}\cr
\displaystyle\prod_{\ell=1}^k \bigl(i - (j-k-1+
\ell)\bigr)/n^k,
\vspace*{2pt}\cr
    \hspace*{32.5pt} \mbox{for $i > j-1$, $j \geq k+2$.}}
\end{equation}
This is a helpful form for bounding the difference between the entries
of $G$ and $H$, which is needed for our convergence analysis in the
next section. Moreover, the functions defined in (\ref{eqhbasis})
give rise to a natural continuous-time parameterization for trend
filtering.

%
\begin{lemma}
\label{lemtfcont}
For inputs $x_1<\cdots<x_n$, and the functions $h_1,\ldots, h_n$
in (\ref{eqhbasis}), define the linear subspace
of functions
%
\begin{equation}
\label{eqhk} \mathcal{H}_k = \operatorname{span}\{h_1,\ldots,
h_n\} = \Biggl\{ \sum_{j=1}^n
\alpha_j h_j\dvtx  \alpha_1,\ldots,
\alpha_n \in \mathbb{R} \Biggr\}.
\end{equation}
If the inputs are evenly spaced, $x_i=i/n$, $i=1,\ldots, n$, then
the continuous-time minimization problem
%
\begin{equation}
\label{eqtfcont} \hat{f}= \argmin_{f \in\mathcal{H}_k} \frac{1}{2}\sum
_{i=1}^n \bigl(y_i-f(x_i)
\bigr)^2 + \lambda\cdot \mathrm{TV}\bigl(f^{(k)}\bigr)
\end{equation}
(where as before, $f^{(k)}$ is understood to mean the $k$th weak
derivative of $f$) is equivalent to the trend filtering problem in
(\ref{eqtf}), that is, their solutions match at the input points,
\[
\hat{\beta}_i = \hat{f}(x_i) \qquad\mbox{for } i=1,
\ldots, n.
\]
\end{lemma}

\begin{figure}

\includegraphics{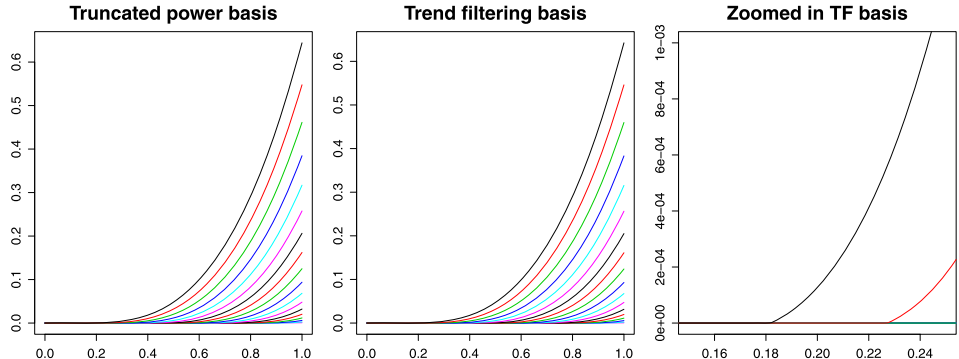}

\caption{For $n=22$ inputs (evenly spaced over
$[0,1]$) and $k=3$, the left panel shows the truncated power basis
functions in \textup{(\protect\ref{eqgbasis})} and the center panel shows the basis
functions in \textup{(\protect\ref{eqhbasis})} utilized by trend
filtering. The two sets of basis functions appear very similar. The
right plot is a zoomed in version of the center plot, and shows the
nonsmooth nature of the trend filtering basis functions---here
(for $k=3$) they have discontinuous first and second derivatives.}\label{figghbasis}
\end{figure}

This result follows by expressing $f$ in (\ref{eqtfcont}) in
finite-dimensional form as $f=\sum_{j=1}^n \alpha_j h_j$, and then
applying\vspace*{-1pt} Lemmas \ref{lemtfh} and \ref{lemtflasso}.
Lemma \ref{lemtfcont} says that the components of trend filtering
estimate, $\hat{\beta}_1,\ldots,\hat{\beta}_n$, can be seen as the
evaluations of a function $\hat{f}\in\mathcal{H}_k$ over the input
points, where
$\hat{f}$ solves the continuous-time problem (\ref{eqtfcont}).
The function $\hat{f}$ is a piecewise polynomial of degree $k$,
with knots contained in $\{x_{k+1},\ldots, x_{n-1}\}$, and for
$k \geq1$, it is continuous since each of the basis functions
$h_1,\ldots, h_n$ are continuous. Hence, for $k=0$ or 1,
the continuous-time trend filtering estimate $\hat{f}$ is a spline (and
further, it is equal to the locally adaptive regression spline
estimate by Lemma \ref{lemsamediff}). But $\hat{f}$ is
not necessarily a spline when $k\geq2$, because in this case it can
have discontinuities in its lower order derivatives (of orders
$1,\ldots, k-1$) at the input
points. This is because each basis function $h_j$, $j =k+2,\ldots, n$,
though infinitely (strongly) differentiable in between the inputs, has
discontinuous derivatives of all lower orders $1,\ldots, k-1$ at the
input point $x_{j-1}$. These discontinuities are visually quite small
in magnitude, and the basis functions $h_1,\ldots, h_n$ look extremely
similar to the truncated power basis functions $g_1,\ldots, g_n$; see
Figure~\ref{figghbasis} for an example.

%

Loosely speaking, the basis functions $h_1,\ldots, h_n$ in
(\ref{eqhbasis}) can be thought of as the falling factorial analogues
of the truncated power basis $g_1,\ldots, g_n$ in (\ref{eqgbasis}).
One might expect then that the subspaces of $k$th degree piecewise
polynomial functions $\mathcal{H}_k$ and $\mathcal{G}_k$ are fairly
close, and that
the (continuous-time) trend filtering and locally adaptive regression
spline problems (\ref{eqtfcont}) and (\ref{eqlrs1}) admit similar
solutions. In the next section, we prove that (asymptotically) this is
indeed the case, though we do so by instead studying the discrete-time
parameterizations of these problems. We show that over a broad
class of true functions $f_0$, trend filtering estimates inherit
the minimax convergence rate of locally adaptive regression splines,
because the two estimators converge to each other at this same
rate.

\section{\texorpdfstring{Rates of convergence.}{Rates of convergence}}\label{secrates}

In this section, we assume an observation model
%
\begin{equation}
\label{eqymodel} y_i = f_0(x_i) +
\varepsilon_i,\qquad i=1,\ldots, n,
\end{equation}
where $x_i=i/n$, $i=1,\ldots, n$ are evenly spaced input points,
$f_0\dvtx [0,1]\rightarrow\mathbb{R}$ is an unknown regression function to
be estimated,
and $\varepsilon_i$, $i=1,\ldots, n$ are i.i.d. sub-Gaussian errors with
zero mean, that is,
%
\begin{eqnarray}\label{eqsubgauss}
\mathbb{E}[\varepsilon_i]=0,\qquad \mathbb{P}\bigl(|
\varepsilon_i|>t\bigr) \leq M \exp \bigl(-t^2/\bigl(2
\sigma^2\bigr) \bigr)
\nonumber\\[-8pt]\\[-8pt]
\eqntext{\mbox{for all } t>0, i=1,\ldots, n}
\end{eqnarray}
for some constants $M,\sigma>0$. [We will write $A_n=O_\mathbb
{P}(B_n)$ to
denote that $A_n/B_n$ is bounded in probability, for random sequences
$A_n,B_n$. We will also write $a_n=\Omega(b_n)$ to denote
$1/a_n=O(1/b_n)$ for constant sequences $a_n,b_n$, and finally
$a_n=\Theta(b_n)$ to denote $a_n=O(b_n)$ and $a_n=\Omega(b_n)$.]

In \citet{locadapt}, the authors consider the same setup,
and study the performance of the locally adaptive regression spline estimate
(\ref{eqlrs1}) when the true function $f_0$ belongs to the set
\[
\mathcal{F}_k(C) = \bigl\{ f\dvtx [0,1]\rightarrow\mathbb{R}\dvtx  \mbox {$f$
is $k$ times weakly differentiable and } \mathrm{TV}\bigl(f^{(k)}\bigr)
\leq C \bigr\}
\]
for some order $k\geq0$ and constant $C>0$. Theorem 10 of
\citet{locadapt} shows that
the $k$th order locally adaptive regression spline estimate $\hat{f}$ in
(\ref{eqlrs1}), with $\lambda=\Theta(n^{1/(2k+3)})$, satisfies
%
\begin{equation}
\label{eqlrsrate} \frac{1}{n}\sum_{i=1}^n
\bigl(\hat{f}(x_i)-f_0(x_i)
\bigr)^2 = O_\mathbb{P}\bigl(n^{-(2k+2)/(2k+3)}\bigr)
\end{equation}
and also that $\mathrm{TV}(\hat{f}) = O_\mathbb{P}(1)$.

\subsection{\texorpdfstring{Minimax convergence rate.}{Minimax convergence rate}}\label{secminimax}

We note that the rate $n^{-(2k+2)/(2k+3)}$ in (\ref{eqlrsrate})
is the minimax rate for estimation over the function class
$\mathcal{F}_k(C)$, provided that $C>1$. To see this, define the Sobolev
smoothness class
\[
\mathcal{W}_k(C) = \biggl\{ f\dvtx [0,1]\rightarrow\mathbb{R}\dvtx  \mbox
{$f$ is $k$ times differentiable and } \int_0^1
\bigl(f^{(k)}(t) \bigr)^2 \,dt \leq C \biggr\}.
\]
Minimax rates over the Sobolev classes are
well-studied, and it is known [e.g., see
\citet{nussbaum}] that
\[
\min_{\hat{f}} \max_{f_0 \in\mathcal{W}_k(C)} \mathbb{E} \Biggl[
\frac{1}{n}\sum_{i=1}^n \bigl(\hat
{f}(x_i)-f_0(x_i) \bigr)^2
\Biggr] = \Omega\bigl(n^{-2k/(2k+1)}\bigr).
\]
Recalling that $\mathrm{TV}(f)=\int_0^1 |f'(t)| \,dt$ for
differentiable $f$, it follows
that\break  $\mathcal{F}_k(C) \supseteq\mathcal{W}_{k+1}(C-1)$, and
\[
\min_{\hat{f}} \max_{f_0 \in\mathcal{F}_k(C)} \mathbb{E} \Biggl[
\frac{1}{n}\sum_{i=1}^n \bigl(\hat
{f}(x_i)-f_0(x_i) \bigr)^2
\Biggr] = \Omega\bigl(n^{-(2k+2)/(2k+3)}\bigr).
\]
In words, one cannot hope to do better than $n^{-(2k+2)/(2k+3)}$ for
function estimation over $\mathcal{F}_k(C)$.

On the other hand, the work of \citet{minimaxwave} provides a
lower bound on the rate of convergence over $\mathcal{F}_k(C)$ for any
estimate linear in $y$---by this, we mean that the vector of its
fitted values over the inputs is a linear function of $y$.
This covers smoothing splines [recall the expression
(\ref{eqssfit}) for the smoothing splines fitted values] and also,
for example, kernel regression estimators. Letting $B^\alpha_{p,q}$ denote
the three parameter Besov space as in \citet{minimaxwave}, and
$\|\cdot\|_{B^\alpha_{p,q}}$ denote the corresponding norm, we have
%
\begin{eqnarray}\label{eqbesovball}
\mathcal{F}_k(C) &\supseteq& \bigl\{ f\dvtx  [0,1] \rightarrow
\mathbb{R}\dvtx  \bigl\|f^{(k)}\bigr\| _\infty+ \mathrm{TV}
\bigl(f^{(k)}\bigr) \leq C \bigr\}\nonumber
\\
&\supseteq& \bigl\{ f\dvtx  [0,1] \rightarrow\mathbb{R}\dvtx
\bigl\|f^{(k)}\bigr\| _{B^1_{1,1}} \leq C' \bigr\}
\\
&\supseteq& \bigl\{ f\dvtx [0,1] \rightarrow\mathbb{R}\dvtx  \|f\|
_{B^{k+1}_{1,1}} \leq C'' \bigr\},\nonumber
\end{eqnarray}
where we write $\|f\|_\infty= \max_{t \in[0,1]} |f(t)|$ for the
$L_\infty$
function norm, and $C',C''$ are constants. The second containment
above follows from a well-known embedding of function spaces [e.g.,
see \citet{wavetour}, \citet{gaussesti}]. The third
\mbox{containment} is given
by applying the Johnen--Scherer bound on the modulus of continuity
[e.g., Theorem 3.1 of \citet{constapprox}] when working with the
usual definition of the Besov norms.\footnote{Thanks to Iain Johnstone
for pointing this out.} Since the minimax linear risk
over the Besov ball in (\ref{eqbesovball})
is of order\break  $n^{-(2k+1)/(2k+2)}$
[\citet{minimaxwave}], we have\footnote{These authors actually
study minimax rates under
the $L_2$ function norm, instead of the discrete (input-averaged) norm
that we consider here. However, these two norms are close enough over
the Besov spaces that the rates do not change; see Section~15.5 of
\citet{gaussesti}.}
\[
\min_{\hat{f}\ \mathrm{linear}} \max_{f_0 \in
\mathcal{F}_k(C)} \mathbb{E} \Biggl[
\frac{1}{n}\sum_{i=1}^n \bigl(\hat
{f}(x_i)-f_0(x_i) \bigr)^2
\Biggr] = \Omega\bigl(n^{-(2k+1)/(2k+2)}\bigr).
\]
Hence, in terms of their convergence rate over $\mathcal{F}_k(C)$, smoothing
splines are suboptimal.

\subsection{\texorpdfstring{Trend filtering convergence rate.}{Trend filtering convergence rate}}\label{sectfrate}

Here, we show that trend filtering also achieves the minimax
convergence rate
over $\mathcal{F}_k(C)$. The arguments used by \citet{locadapt}
for locally
adaptive regression splines cannot be directly applied here, as they
are based
some well-known interpolating properties of splines that do not easily extend
to the trend filtering setting.
Our strategy is hence to show that, as $n\rightarrow\infty$, trend
filtering estimates
lie close enough to locally adaptive regression spline estimates to
share their
favorable asymptotic properties.
(Note that for $k=0$ or $k=1$, the trend filtering and locally adaptive
regression
spline estimates are exactly the same for any given problem instance,
as shown in
Lemma \ref{lemsamediff} in Section~\ref{seclrs}; therefore, the
arguments here are
really directed toward establishing a convergence rate for trend
filtering in the
case $k\geq2$.) Using the triangle inequality (actually, using
$\|a+b\|_2^2 = \|a\|_2^2+\|b\|_2^2 + 2a^Tb \leq2\|a\|_2^2+2\|b\|
_2^2$), we have
%
\begin{equation}
\label{eqtrineq} \qquad\frac{1}{n}\sum_{i=1}^n
\bigl(\hat{\beta}_i-f_0(x_i)
\bigr)^2 \leq \frac{2}{n}\sum_{i=1}^n
\bigl(\hat{\beta}_i-\hat{f}(x_i) \bigr)^2 +
\frac{2}{n}\sum_{i=1}^n \bigl(
\hat{f}(x_i)-f_0(x_i) \bigr)^2,
\end{equation}
where $\hat{\beta}, \hat{f}$ are the trend filtering and locally adaptive
regression spline estimates in (\ref{eqtf}), (\ref{eqlrs1}), respectively.
The second term above is $O_\mathbb{P}(n^{-(2k+2)/(2k+3)})$ by (\ref
{eqlrsrate});
if we could show that the first term above is also $O_\mathbb
{P}(n^{-(2k+2)/(2k+3)})$,
then it would follow that trend filtering converges (in probability) to $f_0$
at the minimax rate.

Recall from Section~\ref{seclrs} that both the trend filtering and locally
adaptive regression spline estimates can be expressed in terms of the fitted
values of lasso problems,
\[
\hat{\beta}= H\hat{\alpha},\qquad \bigl(\hat{f}(x_1),\ldots,
\hat{f}(x_n) \bigr) = G\hat{\theta},
\]
where $G,H \in\mathbb{R}^{n\times n}$ are the basis matrices in (\ref
{eqgg}),
(\ref{eqhh}), and $\hat{\alpha},\hat{\theta}$ are the solutions
in lasso problems
(\ref{eqlassog}), (\ref{eqlassoh}). Hence, we seek a bound for
$\sum_{i=1}^n (\hat{\beta}_i-\hat{f}(x_i))^2 = \|H\hat{\alpha
}-G\hat{\theta}\|_2^2$,
the (squared norm) difference in fitted values between two lasso
problems with the same outcome $y$, but different predictor matrices $G,H$.
Intuitively, a tight bound is plausible here
because $G$ and $H$ have such similar entries (again, for $k=0$ or $k=1$,
we know that indeed $G=H$).

While there are existing results on the stability of the lasso
fit as a function of the outcome vector $y$
[e.g., \citet{lassodf2} show that for any fixed predictor matrix and
tuning parameter value, the lasso fit is nonexpansive as a function of $y$],
as far as we can tell, general stability results do not exist for the
lasso fit as a function of its predictor matrix. To this end, in the
supplement [\citet{trendfilter-supp}], we develop bounds for the
difference in fitted values of two lasso problems that have different
predictor matrices, but the same outcome.
The bounds are asymptotic in nature, and are driven primarily by the
maximum elementwise difference between the predictor matrices. We can
apply this work in the current setting to show that the trend
filtering and locally adaptive\vspace*{1.5pt} regression spline estimates converge
(to each other) at the desired rate,
$n^{-(2k+2)/(2k+3)}$; essentially, this amounts to showing that the
elements of $G-H$ converge to zero quickly enough.

%
\begin{theorem}
\label{thmtflrs}
Assume that $y\in\mathbb{R}^n$ is drawn from the model (\ref{eqymodel}),
with evenly spaced inputs $x_i=i/n$, $i=1,\ldots, n$ and i.i.d. sub-Gaussian errors
(\ref{eqsubgauss}). Assume also that $f_0 \in\mathcal{F}_k(C)$,
that is, for a fixed
integer $k\geq0$ and constant $C>0$, the true function $f_0$ is $k$ times
weakly differentiable and $\mathrm{TV}(f_0^{(k)}) \leq C$. Let $\hat
{f}$ denote the $k$th
order locally adaptive regression spline estimate in (\ref{eqlrs1}) with
tuning parameter $\lambda=\Theta(n^{1/(2k+3)})$, and let $\hat{\beta
}$ denote the
$k$th order trend filtering estimate in (\ref{eqtf}) with tuning parameter
$(1+\delta)\lambda$, for any fixed $\delta>0$. Then
\[
\frac{1}{n}\sum_{i=1}^n \bigl(\hat{
\beta}_i-\hat{f}(x_i) \bigr)^2 =
O_\mathbb{P}\bigl(n^{-(2k+2)/(2k+3)}\bigr).
\]
\end{theorem}

\begin{pf}
We use Corollary 2 in the supplementary document [\citet{trendfilter-supp}], with
$X=G$ and $Z=H$. First, note that our sub-Gaussian assumption in~(\ref{eqsubgauss}) implies that $\mathbb{E}[\varepsilon_i^4]<\infty$
(indeed, it
implies finite moments of all orders), and with $\mu=(f_0(x_1),\ldots,
f_0(x_n))$, we know from the result of \citet{locadapt}, paraphrased
in (\ref{eqlrsrate}), that
\[
\|\mu-G\hat{\theta}\|_2 = O_\mathbb{P}\bigl(n^{-(k+1)/(2k+3)+1/2}
\bigr) = O_\mathbb{P}(\sqrt{n}).
\]
Furthermore, the locally adaptive regression spline estimate $\hat{f}$
has total variation
\[
\mathrm{TV}(\hat{f}) = \|\hat{\theta}_2\|_1 =
O_\mathbb{P}(1),
\]
where $\hat{\theta}_2$ denotes the last $p_2=n-k-1$ components of
$\hat{\theta}$.
Therefore, recalling that $\lambda=\Theta(n^{1/(2k+3)})$, the remaining
conditions needed for Corollary 2 in the supplement [\citet{trendfilter-supp}]
reduce to
\[
n^{(2k+2)/(2k+3)} \|G_2-H_2\|_\infty\rightarrow0
\qquad\mbox {as } n\rightarrow\infty,
\]
where $G_2$ and $H_2$ denote the last $n-k-1$ columns of $G$ and $H$,
respectively, and $\|A\|_\infty$ denotes the maximum absolute element
of a matrix $A$. The above limit can be established by using
Stirling's formula (and controlling the approximation errors) to bound
the elementwise differences in $G_2$ and $H_2$; see Lemma 5
in the supplementary document [\citet{trendfilter-supp}]. Therefore, we apply Corollary
2 in the supplement to conclude that
\[
\|G\hat{\theta}-H\hat{\alpha}\|_2 = O_\mathbb{P}\bigl(\sqrt{n^{1/(2k+3)}}\bigr).
\]
Squaring both sides and dividing by $n$ gives the result.
\end{pf}


Now, using the triangle inequality (\ref{eqtrineq}) [and recalling the
convergence rate of the locally adaptive regression spline estimate
(\ref{eqlrsrate})], we arrive at the following result.

%
\begin{corollary}
\label{cortfrate}
Under the assumptions of Theorem \ref{thmtflrs}, for a tuning parameter
value $\lambda=\Theta(n^{1/(2k+3)})$, the $k$th order trend filtering estimate
$\hat{\beta}$ in (\ref{eqtf}) satisfies
\[
\frac{1}{n}\sum_{i=1}^n \bigl(\hat{
\beta}_i-f_0(x_i) \bigr)^2 =
O_\mathbb{P}\bigl(n^{-(2k+2)/(2k+3)}\bigr).
\]
Hence, the trend filtering estimate converges in probability to $f_0$
at the minimax rate.
\end{corollary}

\begin{rem*}
\citet{locadapt} prove the analogous convergence
result (\ref{eqlrsrate}) for locally adaptive regression splines using
an elegant argument involving metric entropy and the interpolating properties
of splines. In particular, a key step in their proof uses the fact that
for every $k\geq0$, and every function $f\dvtx  [0,1] \rightarrow\mathbb
{R}$ that
has $k$ weak derivatives, there exists a spline $g\in\mathcal{G}_k$
[i.e., $g$ is a~spline of degree $k$ with knots in the set $T$, as
defined in (\ref{eqgk})] such that
%
\begin{equation}
\label{eqgprop} \qquad \max_{x \in[x_1,x_n]} \bigl|f(x)-g(x)\bigr| \leq d_k
\mathrm{TV}\bigl(f^{(k)}\bigr) n^{-k} \quad\mbox{and}\quad
\mathrm{TV}\bigl(g^{(k)}\bigr) \leq d_k \mathrm{TV}
\bigl(f^{(k)}\bigr),
\end{equation}
where $d_k$ is a constant depending only on $k$ (not on the function $f$).
Following this line of argument for trend filtering would require us to
establish the same interpolating properties (\ref{eqgprop}) with
$h\in\mathcal{H}_k$ in place of $g$, where $\mathcal{H}_k$, as
defined in
(\ref{eqhbasis}), (\ref{eqhk}), is
the domain of the continuous-time trend filtering minimization problem
in (\ref{eqtfcont}). This gets very complicated, as $\mathcal{H}_k$
does not
contain spline functions, but instead functions that can have
discontinuous lower order derivatives at the input points $x_1,\ldots, x_n$.
We circumvented such a complication by proving that trend
filtering estimates converge to locally adaptive regression spline
estimates at a rate equal to the minimax convergence rate (Theorem
\ref{thmtflrs}), therefore,
``piggybacking'' on the locally adaptive regression splines rate
due to \citet{locadapt}.
\end{rem*}


\subsection{\texorpdfstring{Functions with growing total variation.}{Functions with growing total variation}}\label{secgrow}

We consider an extension to estimation over the function
class $\mathcal{F}_k(C_n)$, where now $C_n>0$ is not necessarily a
constant and can
grow with $n$. As in the last section, we rely on a result of
\citet{locadapt} for locally adaptive regression splines in the same
situation, and prove that trend filtering estimates and locally
adaptive regression spline estimates are asymptotically very close.

%
\begin{theorem}
\label{thmtflrs2}
Assume that $y\in\mathbb{R}^n$ is drawn from the model (\ref{eqymodel}),
with inputs $x_i=i/n$, $i=1,\ldots, n$ and i.i.d. sub-Gaussian errors
(\ref{eqsubgauss}). Assume also that $f_0 \in\mathcal{F}_k(C_n)$,
that is, for a fixed
integer $k\geq0$ and $C_n>0$ (depending\vspace*{1pt} on $n$), the true function
$f_0$ is
$k$ times weakly differentiable and $\mathrm{TV}(f_0^{(k)}) \leq C_n$.
Let $\hat{f}$ denote
the $k$th order locally adaptive regression\vspace*{1pt} spline estimate in
(\ref{eqlrs1}) with tuning parameter
$\lambda=\Theta(n^{1/(2k+3)}C_n^{-(2k+1)/(2k+3)})$, and let $\hat
{\beta}$
denote the $k$th order trend filtering estimate in (\ref{eqtf}) with
tuning parameter $(1+\delta)\lambda$, for any fixed $\delta>0$.
If $C_n$ does not grow too quickly,
%
\begin{equation}
\label{eqcnsize} C_n=O\bigl(n^{(k+2)/(2k+2)}\bigr),
\end{equation}
then
\[
\frac{1}{n}\sum_{i=1}^n \bigl(\hat{
\beta}_i-\hat{f}(x_i) \bigr)^2 =
O_\mathbb{P}\bigl(n^{-(2k+2)/(2k+3)}C_n^{2/(2k+3)}\bigr).
\]
\end{theorem}

\begin{pf}
The arguments here are similar to the proof of Theorem \ref{thmtflrs}.
We invoke Theorem 10 of \citet{locadapt}, for the present
case of growing total variation $C_n$: this says that
%
\begin{equation}
\label{eqlrsrate2} \frac{1}{n}\sum_{i=1}^n
\bigl(\hat{f}(x_i)-f_0(x_i)
\bigr)^2 = O_\mathbb{P}\bigl(n^{-(2k+2)/(2k+3)}C_n^{2/(2k+3)}
\bigr)
\end{equation}
and also $\mathrm{TV}(\hat{f})=\|\hat{\theta}_2\|_1=O_\mathbb
{P}(C_n)$. Now the conditions for
Corollary 2 in the supplementary document [\citet{trendfilter-supp}] reduce to
%
\begin{eqnarray}
\label{eqcond1}  n^{(2k+2)/(4k+6)}C_n^{(2k+2)/(2k+3)}
\|G_2-H_2\|_\infty &=& O(1),
\\
\label{eqcond2}  n^{(2k+2)/(2k+3)}\|G_2-H_2
\|_\infty&\rightarrow& 0 \qquad\mbox {as } n \rightarrow \infty.
\end{eqnarray}
Applying the assumption (\ref{eqcnsize}) on $C_n$, it is seen that
both (\ref{eqcond1}), (\ref{eqcond2}) are implied by the condition
$n \|G_2-H_2\|_\infty= O(1)$, which is shown in Lemma 5
in the supplementary document [\citet{trendfilter-supp}].
Therefore, we conclude using Corollary 2 in the
supplement that
\[
\sqrt{\sum_{i=1}^n \bigl(\hat{
\beta}_i-\hat{f}(x_i) \bigr)^2} = \|H\hat{
\alpha}-G\hat{\theta}\|_2 = O_\mathbb{P} \Bigl(\sqrt {n^{1/(2k+3)}C_n^{2/(2k+3)}}
\Bigr),
\]
which gives the rate in the theorem after squaring both sides and
dividing by $n$.
\end{pf}

Finally, we employ the same triangle inequality (\ref{eqtrineq}) [and
the locally adaptive regression splines result (\ref{eqlrsrate2}) of
\citet{locadapt}] to derive a rate for trend filtering.

%
\begin{corollary}
\label{cortfrate2}
Under the assumptions of Theorem \ref{thmtflrs2}, for
$C_n=\break  O(n^{(k+2)/(2k+2)})$ and a tuning parameter value
$\lambda=\Theta(n^{1/(2k+3)}C_n^{-(2k+1)/(2k+3)})$,
the $k$th order trend filtering estimate $\hat{\beta}$ in
(\ref{eqtf}) satisfies
\[
\frac{1}{n}\sum_{i=1}^n \bigl(\hat{
\beta}_i-f_0(x_i) \bigr)^2 =
O_\mathbb{P}\bigl(n^{-(2k+2)/(2k+3)}C_n^{2/(2k+3)}\bigr).
\]
\end{corollary}

\begin{rem*}
Although we manage to show that trend filtering achieves
the same convergence rate as locally adaptive regression splines in
the case of underlying functions with growing total variation, we require
the assumption that $C_n$ grows no faster than $O(n^{(k+2)/(2k+2)})$, which
is not required for the locally adaptive regression spline result
proved in
\citet{locadapt}. But it is worth pointing out that for $k=0$ or $k=1$,
the restriction $C_n=O(n^{(k+2)/(2k+2)})$ for the trend filtering convergence
result is not needed, because in these cases trend filtering and
locally adaptive
regression splines are exactly the same by Lemma \ref{lemsamediff} in
Section~\ref{seclrs}.
\end{rem*}

\section{\texorpdfstring{Astrophysics data example.}{Astrophysics data example}}\label{secastro}

\begin{figure}

\includegraphics{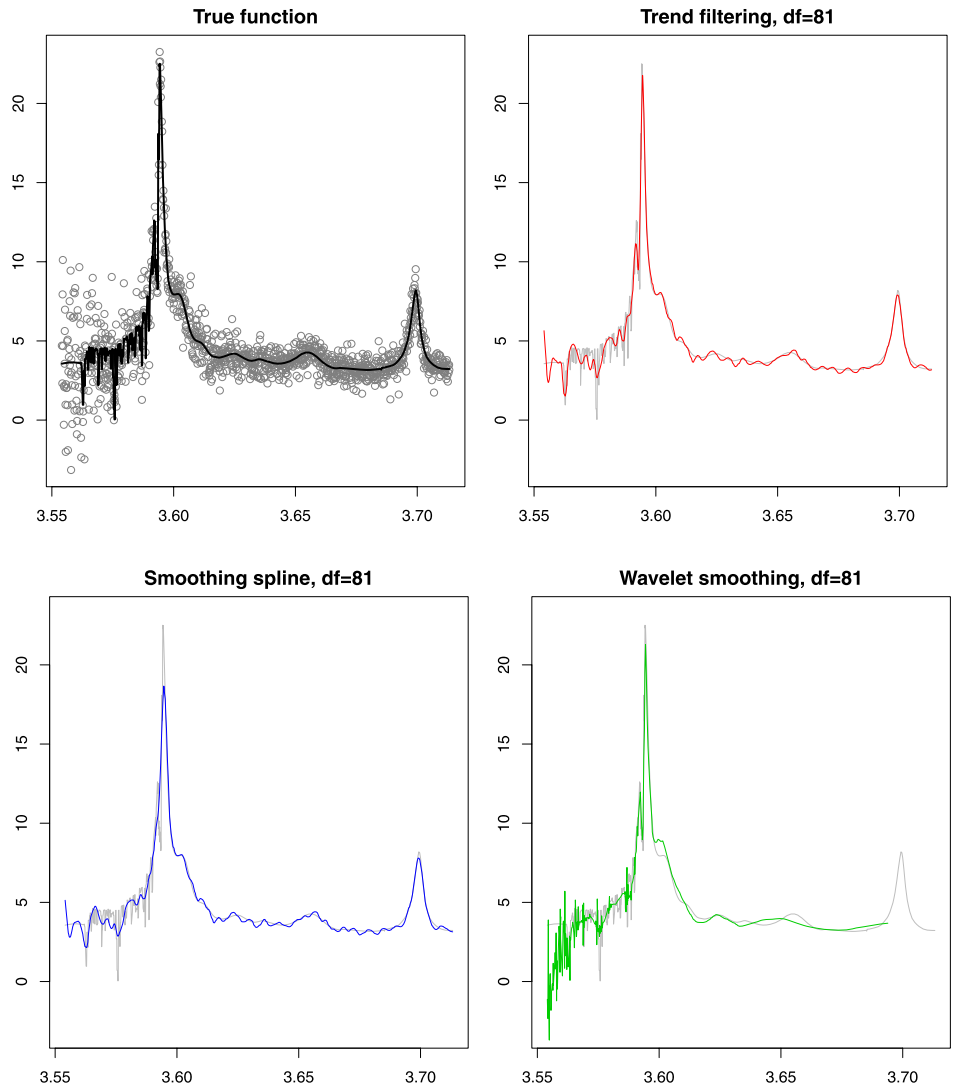}

\caption{The top left panel shows data simulated from a model
for quasar spectrum. The true curve, in black, displays flux as a
function of log wavelength. The gray points are noisy
observations at $n=1172$ wavelength values. We fit trend
filtering, smoothing splines, and wavelets to these points,
and tuned each to have 81 degrees of freedom. (This value was
chosen because it corresponded to the trend
filtering model with the minimum squared error loss, averaged 20
simulations---see Figure~\protect\ref{figspectrumerr}.) The resulting
estimates are displayed in the top right, bottom
left and bottom right panels, respectively (with the true
function plotted in the background, in gray). Trend filtering and
smoothing splines give similar fits, except that trend filtering
does a better job of estimating the large peak at around $x=3.6$, as
well as some of the finer features of the true function to the left
of this. Wavelet
smoothing also does well in detecting the extreme peak, but then
overfits the true function on the left-hand side, and underfits on the
right.}\label{figspectrum}
\end{figure}

We examine data from an astrophysics simulation model for quasar
spectra, provided by Yu Feng, with help from Mattia Ciollaro and
Jessi Cisewski.
Quasars are among the most luminous objects in the universe.
Because of this, we can observe them at great distances, and features
in their spectra reveal information about the universe along the
line-of-sight of a given quasar. A quasar spectrum drawn from this
model is displayed in the top left panel of Figure~\ref{figspectrum}.
The spectrum, in black, shows the flux ($y$-axis) as
a function of log wavelength ($x$-axis). Noisy realizations of this
true curve are plotted as gray points, measured
at $n = 1172$ points, (approximately) equally spaced on the log
wavelength scale. We see that the true function has a
dramatically\vadjust{\goodbreak}
different level of smoothness as it traverses the wavelength scale,
being exceptionally wiggly on the left-hand side of the domain, and much
smoother on the right.
(The error variance is also seen to be itself inhomogeneous, with larger
errors around the wiggly portions of the curve, but for simplicity we do
not account for this.) The wiggly left-hand side of the spectrum are
absorption features called the Lyman-alpha forest.

To estimate the underlying function, we applied trend filtering,
smoothing splines, and wavelet smoothing, each over 146 values of
degrees of freedom (from~4 to 150 degrees to freedom). Locally
adaptive regression splines were not compared because of their
extreme proximity to trend filtering. The \texttt{smooth.spline}
function in R was used to fit the smoothing spline estimates, and
because it produces cubic order smoothing splines, we considered
cubic order trend filtering and wavelets with 4 vanishing moments
to put all of the methods on more or less equal footing. Wavelet
smoothing was fit using the \texttt{wavethresh} package in R, and
the ``wavelets on the interval'' option was chosen to handle
the boundary conditions (as periodicity and symmetry are not
appropriate assumptions for the boundary behavior in this example),
which uses an algorithm of \citet{waveinterval}. Wavelet transforms
generally require the number of observations to be a power of 2
(this is at least true in the \texttt{wavethresh} implementation),
and so we restricted the wavelet smoothing estimate to use the first
1024 points with the smallest log wavelengths.

Figure~\ref{figspectrum} demonstrates the function estimates from
these methods, run on the single data set shown in the top left
panel (the observations are not drawn in the remaining panels so as
not to cloud the plots). Each estimate was tuned to have 81 degrees
of freedom. We can see that trend filtering (top right panel)
captures many features of the true function, picking up the large
spike just before $x=3.6$, but missing some of the action on the
left-hand side. The smoothing spline estimate (bottom left) appears
fairly similar, but it does not fit the magnitudes of the
wiggly components as well. Wavelet
smoothing (bottom right) detects the large spike, but badly overfits the
true function to the left of this spike, and even misses gross smoothness
features to the right.

%

\begin{figure}

\includegraphics{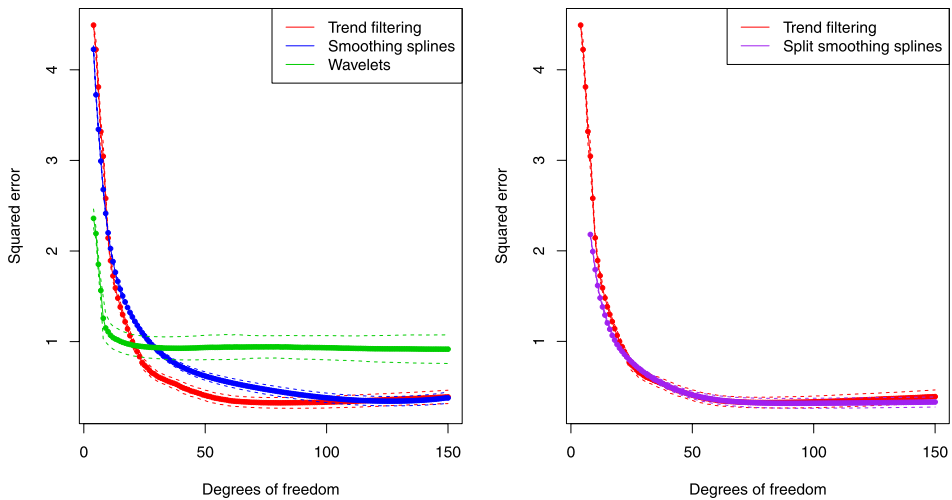}

\caption{Left plot: the squared error
loss in estimating the true light curve for the quasar spectrum data,
using trend filtering, smoothing splines and wavelets, fit over a
range of model complexities (degrees of freedom values).
The results were averaged over 20 simulated data sets, with standard
errors drawn as dotted lines. Trend filtering achieves a significantly
lower squared error loss than smoothing splines for models of low
complexity; both perform much better than wavelets. Right plot: trend
filtering versus a smoothing spline estimator that fits a
different smoothing parameter on two halves of the domains (on either
side of the peak around $x=3.6$); the methods perform
comparably.}
\label{figspectrumerr}
\end{figure}

We further compared the three contending methods by computing
their average squared error loss to the true function, over 20
draws from the simulated model. This is shown in the left panel
of Figure~\ref{figspectrumerr}. Trend filtering outperforms
smoothing splines for lower values of model complexity (degrees
of freedom); this can be attributed to its superior capability
for local adaptivity, a claim both empirically supported by the
simulations in Section~\ref{secssempir}, and formally explained
by the theory in Section~\ref{secrates}. Wavelet smoothing is
not competitive in terms of squared error loss. Although
in theory it
achieves the same (minimax) rate of convergence as trend filtering,
it seems in the current setting to be hurt by the
high noise level at the left-hand side of the domain; wavelet smoothing
overfits in this region, which inflates the estimation variance.

%

Finally, we compared trend filtering to a smoothing spline estimator
whose tuning parameter varies over the input domain
(to yield a finer level of local adaptivity). For an example of a recent
proposal of such an estimator, see \citet{smoothsplinesvary}
(see also the references therein). Methods
that fit a flexibly varying tuning parameter over the domain can
become very complicated, and so to simplify matters for the quasar
spectrum data, we allowed the smoothing spline estimator two different
tuning parameters $\lambda_1,\lambda_2$ to the left and right of
$x=3.6$. Note that this represents somewhat of an ideal scenario
for variable parameter smoothing splines, as we fixed an
appropriate division of the domain based on knowledge of the true
function. It should also be noted that we fit the split smoothing
spline estimator over a total of $146 \cdot146 ={}$21,316 values of
degrees of freedom (146 in each half of the domain), which puts it
at an advantage over the other methods. See
the right panel of Figure~\ref{figspectrumerr} for the results.
Over 20 simulated data sets, we fit split smoothing splines whose
degrees of freedom $d_1,d_2$ on the left and right sides of $x=3.6$
ranged from 4 to 150. For each value of degrees of freedom
$d$, the plotted curve shows the minimum squared error loss over all models
with $d_1+d_2=d$. Interestingly, despite all of the advantages
imparted by its setup, the split smoothing spline estimator performs
basically on par with trend filtering.

\section{\texorpdfstring{Extensions and discussion.}{Extensions and discussion}}\label{secextdisc}

We have shown that trend filtering, a newly proposed method for
nonparametric regression of \citet{l1tf}, is both
fast and locally adaptive. Two of the major tools for adaptive
spline estimation are smoothing splines and locally adaptive
regression splines; in short, the former estimators are fast but not
locally adaptive, and the latter are locally adaptive but not fast.
Trend filtering lies in a comparatively favorable position: its
estimates can be computed in $O(n^{3/2})$ worst-case complexity
(at a fixed value of the tuning parameter $\lambda$, using a
primal--dual interior point algorithm), which is slower than the $O(n)$
complexity of smoothing splines, but not by a big margin; its
estimates also achieve the same convergence rate as locally adaptive
regression splines over a broad class of underlying functions
(which is, in fact, the minimax rate over this class).\looseness=-1

One way to construct trend filtering estimates, conceptually, is to
start with the lasso form for locally adaptive regression
splines (\ref{eqlassog}), but then replace the matrix $G$ in
(\ref{eqgg}), which is generated by the truncated power series,
with the matrix $H$ in~(\ref{eqtfh}), generated by
something like their falling factorial
counterparts. This precisely defines trend filtering, and it
has the distinct computational
advantage that $H$ has a sparse banded inverse (whereas the inverse of
$G$ is dense). Moreover, the matrix $H$ is close enough to $G$ that
trend filtering estimates retain some of the desirable theoretical
properties of locally adaptive regression splines, that is, their minimax
rate of convergence. Although this change-of-basis perspective is
helpful for the purposes of mathematical analysis, the original
representation for trend filtering (\ref{eqtf}) is certainly more
natural, and also perhaps more useful for constructing related
estimators whose characteristics go beyond (piecewise) polynomial
smoothness of a given order. We finish by discussing this, in
Section~\ref{secsva}. First, we briefly discuss an extension to multivariate
inputs.

\subsection{\texorpdfstring{Multivariate trend filtering.}{Multivariate trend filtering}}\label{secmulti}

An important extension concerns the
case of multivariate inputs $x_1,\ldots, x_n \in\mathbb{R}^p$. In this case,
there are two strategies for extending trend filtering that one
might consider. The first is to extend the definition of the discrete
difference operators to cover multivariate inputs---the analogous
extension here for smoothing splines are thin plate splines
[\citet{wahbasplines,greensilver}]. An extension such as this is
``truly'' multivariate, and is an ambitious undertaking; even just the
construction of an appropriate multivariate discrete difference
operator is a topic deserving its own study.

A second, more modest approach for multivariate input points is to
fit an additive model whose individual component functions are fit
by
(univariate) trend filtering.\vadjust{\goodbreak} \citet{gam} introduced additive models,
of the form
%
\begin{equation}
\label{eqymodel2} y_i = \sum_{j=1}^p
f_j(x_{ij}) + \varepsilon_i,\qquad i=1,\ldots, n.
\end{equation}
The model (\ref{eqymodel2}) considers the contributions from each
variable in the input space marginally. Its estimates will
often scale better with the underlying dimension $p$, both in terms of
computational and statistical efficiency, when compared to those from
a ``true'' multivariate extension that considers
variables jointly. Fitting the component functions $\hat{f}_1,\ldots,
\hat{f}_p$ is most often done with a backfitting (blockwise coordinate
descent) procedure, where we cycle through estimating each
$\hat{f}_j$ by fitting the current residual to the $j$th variable,
using a univariate nonparametric regression estimator. Common
practice is to use smoothing splines for these individual univariate
regressions, but given their improved adaptivity properties and
comparable computational efficiency, using trend filtering estimates
for these inner regressions is an idea worth investigating.

\subsection{\texorpdfstring{Synthesis versus analysis.}{Synthesis versus analysis}}\label{secsva}

Synthesis and analysis are concepts from signal processing that,
roughly speaking,
describe the acts of building up an estimator by adding together a
number of fundamental components, respectively, whittling down an
estimator by removing certain undesirable components. The same terms
are also used to convey related concepts in many scientific fields.
In this section, we compare synthesis and analysis in the context of
$\ell_1$ penalized estimation. Suppose that we want to construct an
estimator of $y\in\mathbb{R}^n$ with some particular set of desired
characteristics, and consider the following two general problem
formulations:
%
\begin{eqnarray}
\label{eqsynthesis} & \displaystyle\min_{\theta\in\mathbb{R}^p} \frac{1}{2}\|y-X\theta
\|_2^2 + \lambda\|\theta\|_1,&
\\
& \displaystyle \min_{\beta\in\mathbb{R}^n} \frac{1}{2}\|y-\beta\|_2^2
+ \lambda\|D\beta\|_1.& \label{eqanalysis}
\end{eqnarray}
The first form (\ref{eqsynthesis}) is the synthesis
approach: here, we choose a matrix $X \in\mathbb{R}^{n\times p}$ whose
columns are atoms or building blocks for the characteristics that we
seek, and in solving the synthesis problem (\ref{eqsynthesis}), we
are adaptively selecting a number of these atoms to form our
estimate of $y$. Problem (\ref{eqanalysis}), on the other
hand, is the analysis approach: instead of enumerating
an atom set via $X$, we choose a penalty matrix $D \in\mathbb
{R}^{m\times n}$
whose rows represent uncharacteristic behavior. In solving the problem
(\ref{eqanalysis}), we are essentially orthogonalizing our estimate with
respect to some adaptively chosen rows of $D$, therefore,
directing it away from uncharacteristic behavior.

The original representation of trend filtering in
(\ref{eqtf}) falls into the analysis framework, with $D=D^{(k+1)}$,
the $(k+1)$st order discrete difference operator; its basis
representation in (\ref{eqlassoh}) falls into the synthesis
framework, with \mbox{$X=H$},\vadjust{\goodbreak} the falling factorial basis matrix (an
unimportant difference is that the $\ell_1$ penalty only applies to
part of the coefficient vector). In the former, we shape
the trend filtering estimate by penalizing jumps in its $(k+1)$st
discrete derivative across the input points; in the latter, we
build it from a set of basis functions, each of which is
nonsmooth at only one different input point. Generally, problems
(\ref{eqsynthesis}) and (\ref{eqanalysis}) can be equated if $D$ has
full row rank (as it does with trend filtering), but not
if $D$ is row rank deficient
[see \citet{genlasso}, \citet{avs}].

%
\begin{figure}

\includegraphics{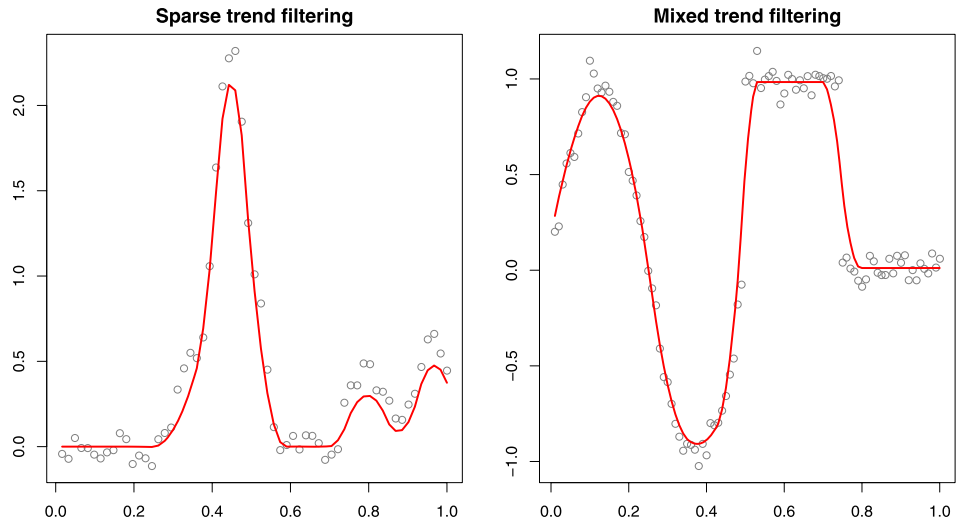}

\caption{\it Left panel: a small example of sparse
quadratic trend filtering $(k=2)$. The estimate $\hat{\beta}$ in~\textup{(\protect\ref{eqstf})} is identically zero for inputs approximately
between $0$
and $0.25$, and $0.6$ and $0.75$. Right panel: an example of
constant/quadratic mixed trend filtering $(k_1=0$ and $k_2=2)$.
The estimate defined in~\textup{(\protect\ref{eqmtf})} is first piecewise quadratic
over the first half of its domain, but then is flat in two stretches
over the second half.}
\label{figanalysis}
\end{figure}

Here, we argue that it can actually be easier to work from the
analysis perspective instead of the synthesis perspective for the
design of nonparametric regression estimators. (The reverse can also
be true in other situations, though that is not our focus.)
For example, suppose that we wanted to construct an estimator that
displays piecewise polynomial smoothness across the input points,
but additionally, is identically zero over some appropriately chosen
subintervals in its domain. It helps to see an example:
see the left panel in Figure~\ref{figanalysis}. Working
from the analysis point of view, such an estimate is easily
achieved by adding a pure $\ell_1$ penalty to the usual trend
filtering criterion, as in
%
\begin{equation}
\label{eqstf} \hat{\beta}= \argmin_{\beta\in\mathbb{R}^n} \frac{1}{2}\|y-\beta
\|_2^2 + \lambda_1 \bigl\|D^{(k+1)}\beta
\bigr\|_1 + \lambda_2 \|\beta\|_1.
\end{equation}
We call (\ref{eqstf}) the sparse trend filtering estimate. This could
be of interest if, for example, $y\in\mathbb{R}^n$ is taken to be the pairwise
differences between two sequences of observations, for example, between
two response curves over time; in this case, the zeros of $\hat{\beta}$
indicate regions in time over which the two responses are deemed to be
more or less the same. It is important to note that an estimate with
these properties seems difficult to construct from the
synthesis perspective---it is unclear what basis
elements, when added together, would generically yield an estimate
like that in~(\ref{eqstf}).

As another example, suppose that we had prior belief that the
observations $y\in\mathbb{R}^n$ were drawn from an underlying
function that
possesses different orders of piecewise polynomial smoothness, $k_1$
and $k_2$, at different parts of its domain. We could then solve the
mixed trend filtering problem,
%
\begin{equation}
\label{eqmtf} \hat{\beta}= \argmin_{\beta\in\mathbb{R}^n} \frac{1}{2}\|y-\beta
\|_2^2 + \lambda_1 \bigl\|D^{(k_1+1)}\beta
\bigr\|_1 + \lambda_2 \bigl\|D^{(k_2+1)}\beta\bigr\|_1.
\end{equation}
The right panel of Figure~\ref{figanalysis} shows an example, with an
underlying function that is mixed piecewise quadratic and piecewise
constant. Again it seems much more difficult to construct an estimate
like (\ref{eqmtf}), that is, one that can flexibly adapt to the
appropriate order of smoothness at different parts of its
domain, using the synthesis framework. Further study
of the synthesis versus analysis perspectives for estimator
construction will be pursued in future work.

\section*{\texorpdfstring{Acknowledgements.}{Acknowledgements}}
We thank the reviewing Editors and referees for their
comments, which lead to an improved paper. We thank Alessandro Rinaldo
and Larry Wasserman for many
helpful and encouraging conversations. We thank Iain Johnstone for
his generous help regarding Besov spaces and the corresponding minimax
rates, in Section~\ref{secminimax}. We also thank
Mattia Ciollaro and Jessi Cisewski for their help with the
astrophysics example in Section~\ref{secastro}. Finally, we thank
Trevor Hastie and Jonathan Taylor for early stimulating conversations
that provided the inspiration to work on this project in the first place.

\begin{supplement}[id=suppA]
\stitle{Supplement to ``Adaptive piecewise polynomial estimation via trend filtering''}
\slink[doi]{10.1214/13-AOS1189SUPP} 
\sdatatype{.pdf}
\sfilename{AOS1189\_supp.pdf}
\sdescription{We provide proofs for the results in Sections~\ref
{seclrs} and \ref{seccont}. We also present the underlying
theoretical framework needed to establish the convergence rates in
Section~\ref{secrates}. 
Finally, we discuss an extension of trend filtering to the case of
arbitrary input points.}
\end{supplement}

%

\printaddresses

\end{document}